\ifx\LocalElevenPtArticleFormatLoaded\undefined
\documentclass[a4paper,11pt]{article} % `article' for A4 paper
\usepackage{array} % array extensions
\usepackage{theorem} % theorem extensions
\usepackage{amsmath,amscd,amssymb} % AMS-LaTeX with CD and symbols
\usepackage{latexsym} % old LaTeX symbols
\usepackage{graphics}
\usepackage{epsfig}
\usepackage{xypic}
\fi

\newtheorem{theorem}{Theorem}[section]
\usepackage{theorem}

\newtheorem{lemma}[theorem]{Lemma}
\newtheorem{corollary}[theorem]{Corollary}
\newtheorem{proposition}[theorem]{Proposition}

\theorembodyfont{\rm}
\newtheorem{remark}{Remark}
\newtheorem{remarks}{Remarks}

\newtheorem{definition}{Definition}
\newenvironment{Proof}{\begin{ProofwCaption}{Proof}}{\end{ProofwCaption}}
\newenvironment{Proof*}[1]{\begin{ProofwCaption}{{#1}}}{\end{ProofwCaption}}
\newenvironment{ProofwCaption}[1]%
{\addvspace\theorempreskipamount \noindent{\it #1.}\rm}%
{\qed \par \addvspace\theorempostskipamount}
\newcommand{\qedsymbol}{\mbox{$\Box$}}
\newcommand{\qed}{\quad\qedsymbol}
\newcommand{\acts}{\curvearrowright}
\newcommand{\wG}{\widetilde{G}}
\newcommand{\oRg}{\overline{\mathcal R}_{g}}
\newcommand{\tu}{{}^t\!}
\newcommand{\Avor}{\overline{A}^{\operatorname{Vor}}}

\newcommand{\Spec}{\operatorname{Spec}}
\newcommand{\Iso}{\operatorname{Iso}}
\newcommand{\Hom}{\operatorname{Hom}}
\newcommand{\id}{\operatorname{id}}
\newcommand{\inv}{^{-1}}

\newcommand{\supp}{\operatorname{supp}}

\newcommand{\GG}{{\mathbb{G}}}

\newcommand{\PP}{{\mathbb{P}}}
\newcommand{\QQ}{{\mathbb{Q}}}
\newcommand{\RR}{{\mathbb{R}}}
\newcommand{\TT}{{\mathbb{T}}}
\newcommand{\ZZ}{{\mathbb{Z}}}
\newcommand{\MM}{{\mathbb{M}}}

\newcommand{\XX}{{\mathbb{X}}}

\newcommand{\calO}{{\cal O}}

\newcommand{\calL}{{\cal L}}
\newcommand{\calM}{{\cal M}}
\newcommand{\calP}{{\cal P}}

\newcommand{\calR}{{\cal R}}

\newcommand{\Cone}{\operatorname{Cone}}
\newcommand{\Proj}{\operatorname{Proj}}
\newcommand{\Pic}{\operatorname{Pic}}
\newcommand{\im}{\operatorname{im}}
\newcommand{\coker}{\operatorname{coker}}
\newcommand{\Nm}{\operatorname{Nm}}

\begin{document}
\begin{center}
{\bf\Large Degenerations of Prym Varieties}\\[3mm]
{\large V. Alexeev, Ch. Birkenhake and K. Hulek}
\end{center}

\noindent {\bfseries Abstract.}  {\footnotesize
Let $(C,\iota)$ be a stable curve with an involution. Following a
classical construction one can define its Prym variety $P$, which in
this case turns out to be a semiabelian group variety and usually not
complete. In this paper we study the question whether there are ``good''
compactifications of $P$ in analogy to compactified Jacobians. The answer to
this question depends on whether we consider degenerations of principally
polarized Prym varieties or degenerations with the induced (non-principal)
polarization. We describe degeneration data of such degenerations. The main
application of our theory lies in the case of degenerations of principally
polarized Prym varieties where we ask whether such a degeneration depends on a
given one-parameter family containing $(C,\iota)$ or not. This allows us to
determine the indeterminacy locus of the Prym map.

\tableofcontents}
\section{Introduction}

Classically, for a smooth projective curve $C$ one defines a
principally pola\-rized abelian variety, its Jacobian $JC$. If $C$ has
nodes, one has, in fact, two analogues: the Picard variety, which in
this case is a semiabelian group variety and usually not complete, and
a "compactified Jacobian", which is a projective variety.  Oda and
Seshadri \cite{OS} define several such compactifications.  Namikawa
\cite{Nam2} defines one variety $\overline{J}C$, and in \cite{A1} it is shown
how to construct a theta divisor $\Theta$ on it and how to obtain the
pair $(\overline{J},\Theta)$ as a stable semiabelic pair.  According
to Mumford and Namikawa \cite{Nam2}, by using methods of toric
geometry, the Torelli map from the moduli space of curves $M_g$ to the
moduli space $A_g$ of principally polarized abelian varieties can be
extended to a morphism from the Deligne-Mumford compactification
$\overline{M}_g$ to the toroidal compactification of $A_g$ for the 2nd
Voronoi fan. The latter compactification
appears in \cite{A2} as the closure of $A_g$
in the moduli of stable semiabelic pairs and \cite{A1} gives a moduli
interpretation of the extended map.

Again, classically, a smooth projective curve $C$ with an involution $\iota$
defines a polarized abelian variety, the Prym variety $P(C,\iota)$.  When the
involution has 0 or 2 fixed points, the polarization is (twice) a
principal polarization. What happens if we now consider a nodal curve
$(C,\iota)$ with an involution?  The quotient curve $C'=C/\iota$ is
again a nodal curve and the identity component of the kernel of the
norm map $\operatorname{Nm}:JC \to JC'$ is a semiabelian variety. One is
naturally led to the question whether this variety has "good"
compactifications.
The naive approach would be to take the closure of the identity component of
the kernel of the norm map in $\overline{J}C$. Examples quickly show that
this naive approach does not always give the right answer. Indeed, the question
allows more than one interpretation. The restriction of the theta divisor of
a Jacobian to the Prym in the case of a fixed point free involution of a
smooth curve defines twice a principal polarization on the Prym variety.
On a smooth abelian variety this defines uniquely a principal polarization.
The situation is more complicated for degenerations of abelian varieties.
There are very well behaved degenerations of polarized abelian varieties
with a divisible polarization such that the limit polarization is not divisible
(e.g. the degeneration of elliptic curves with a degree $n$ divisor to
a fibre of type $I_n$). Hence we must expect different answers whether we
consider the degeneration problem for principally polarized Pryms or for
Prym varieties with the induced polarization. We will address both of these
questions in this paper.

The first question is closely related to the Prym map from the
space $\calR_g$ of
unramified double covers of a curve of genus $g$ to $A_{g-1}$. The space
$\calR_g$ has a natural normal compactification  $\overline{\calR}_g$
parametrizing
admissible double covers of stable curves and the question arises whether the
Prym map extends to a regular morphism from $\overline{\calR}_g$ to the
second Voronoi compactification of $A_{g-1}$ (or some other suitable
compactification). Beauville has constructed a partial extension of the
Prym map in the case where the kernel of the norm map has no toric part, i.e.
where the image of the Prym map is still contained in $A_{g-1}$. The preprint
version of the paper by Friedman and
Smith \cite{FS} contains some examples of admissible double covers where
the Prym map does not extend to any reasonable
toroidal compactification of $A_{g-1}$. The main result of the present
paper is a
combinatorial criterion (see Theorem \ref{thm:indeterminacy_locus})
which characterizes those admissible double covers
which lie in the indeterminacy locus of the Prym map. In the appendix to this
paper by Vologodsky \cite{V} this combinatorial criterion is translated into a
more geometric criterion which says that the locus of indeterminacy is the
closure of the locus given by the examples found by Friedman and Smith.
It should be pointed out, that the fact that the Prym map has indeterminacies,
is not related to singularities of $\overline{\calR}_g$. Indeed, the variety
$\overline{\calR}_g$ is smooth at a generic point representing an
example of Friedman-Smith type.

The contents of this paper is organized as follows. In
Section~1 we study the kernel of the norm map
$\operatorname{Nm}:JC \to JC'$. Its identity component is a
semiabelian variety $P$, whose classifying
homomorphism we determine.
Moreover we
compute the number of components of the kernel of the norm map (Proposition
\ref{prop16}) and determine the induced polarization type of the
abelian part of $P$ (Lemma \ref{lem17}, Proposition
\ref{prop18}, and the remarks following these statements).

Section~2 contains, for the reader's convenience, a brief summary of the
relevant parts of the degeneration theory of abelian varieties in a form
which is relevant for this paper. What is new here is a discussion of induced
degeneration data.

In Section~3 we consider degenerations of Jacobian varieties and of
principally polarized Prym varieties.
Here we shall relate the degeneration data which we obtain to the results
from section~1.
As the main application we determine
the precise locus of indeterminacy of the extended Prym map
(Theorem~\ref{thm:indeterminacy_locus}).
Our answer is given in terms of
combinatorial data of the graph $(\Gamma, \iota)$ associated to the curve $C$.
Given a concrete example
it is easy to check whether the point $[(C,\iota)]$ lies in the
indeterminacy locus or not.

In Section~4 we consider the non-principally polarized case. For any
nodal curve with an involution (including those to which the Prym map
does not extend) we construct a ``middle'' compactified Prym variety of which
we show that it is the limit of non-principally polarized
Pryms (Theorem \ref{theoremmiddlePrym}).
This construction is closer in spirit to the naive approach of
taking the closure of
$P$ in $\overline{J}C$. It is, therefore, natural to compare
this compactification of the Prym variety to the compactified Jacobian
$\overline{J}C$. In order to do this we construct
several finite morphisms $f_v$ from the ``middle'' compactified Prym
variety to the compactified Jacobian of $C$ (see Theorem
\ref{theoremembedding}).  Among these morphisms we identify the best
one, which corresponds to a ``maximal half-integral shift'' $v$.

Section~5 is devoted to a series of examples which illustrate various
aspects of our theory. In particular, example \ref{F-S} is devoted to the
Friedman-Smith examples.
In the simplest case where the Prym map is not defined we obtain a
${\PP}^1$ of different degenerations. In example \ref{3-comp} we illustrate
in a
simple case how the degenerations in the principally polarized and in the
non-principally polarized case are related to each other.

{\bf Acknowledgement.}
We would like to thank Professors R.Smith and R.Varley for helpful discussions.
The first named author was supported by the NSF.
The last named author is grateful to the
University of Georgia at Athens for hospitality during two visits.
This project was further
supported by the DFG special research project
``Global methods in complex geometry'', grant Hu 337/5-1.

\section{The semiabelian part of the Prym variety}
\subsection{Construction of $P$}\label{Construction of $P$}

Throughout this paper (with the exception of sections $2$ and $3$ where we
will work with families) we will denote by $C$ a connected {\em nodal}
curve
defined over an algebraic field $k$ of characteristic different from $2$
, i.e. a projective curve whose singularities are at
most nodes.

To every such curve we can associate a graph
$\Gamma=\Gamma(C)$ whose vertices $\{v_i\}_{i\in I}$ correspond to the
irreducible components $C_i$ of $C$ and whose edges $\{e_j\}_{j\in J}$
correspond to the nodes $Q_j$ of $C$. After choosing an
orientation (it will
not matter which one we choose) for $\Gamma$ one can define a chain complex
$C_{\bullet} (\Gamma,
\ZZ)$ where
$$
C_0(\Gamma,\ZZ)=\bigoplus_{i\in I} \ZZ v_i ,\quad
C_1(\Gamma,\ZZ)=\bigoplus_{j\in J} \ZZ e_j
$$
and where
$
\partial (e_j) = v_i - v_{i'}
$
if $e_j$ is an edge going from $v_i$ to $v_{i'}$. Moreover we choose pairings
$[ , ]$ on $C_0(\Gamma,\ZZ)$ and $( , )$ on $C_1(\Gamma, \ZZ)$ by
$[v_i, v_j]=\delta_{ij}$ and $(e_i, e_j)=\delta_{ij}$. One can then identify
the cochain complex $C^{\bullet}(\Gamma, \ZZ)$ with the chain complex
$C_{\bullet}(\Gamma,\ZZ)$ via these pairings. Thus the coboundary
map $\delta$ becomes
$\delta v_i=\sum\limits_{j\in J} [v_i, \partial e_j] e_j$. The homology and
cohomology groups $H_i(\Gamma,\ZZ)$ and
$H^i(\Gamma,\ZZ); i=0,1$ are defined in the usual way. For details of this
see
\cite[chapter I]{OS}.

The Jacobian $JC$ is defined as the group of line bundles on $C$
whose multidegree is $0$. This defines a group scheme which coincides with
the usual
Jacobian if $C$ is smooth, but which is not necessarily compact otherwise. We
denote by
$\nu: N\rightarrow C$ the normalization of $C$. By \cite
[Proposition 10.2]{OS} the group scheme
$JC$ is an extension
\begin{eqnarray}
1\longrightarrow H^1(\Gamma,k^{\ast})\longrightarrow JC
\stackrel{\nu^{\ast}}{\longrightarrow} JN \longrightarrow 0
\end{eqnarray}
where $JN = \prod\limits_{i\in I}JN_i$ is
the product of the Jacobians of the normaliza\-tions of the components $C_i$ of
$C$. This extension defines an element in the group
$\mbox{Ext}^1(JN, H^1 (\Gamma, k^{\ast }))
= \mbox{Hom} (H_1(\Gamma,\ZZ ), {^tJ}N).$

It is not difficult to make this explicit: Every edge $e_j$ corresponds to a
double point $Q_j$ and over $Q_j$ lie two points $Q_j^{+}$ and $Q_j^{-}$ in
$N$. If $e_j$ is not a loop, i.e. an edge going from $v_i$ to $v_{i'}$ with
$i'\neq i$ we can distinguish $Q_j^{+}$ and $Q_j^{-}$ by saying that $Q_j^{+}$
lies on
$N_i$ and $Q_j^{-}$ lies on $N_{i'}$. Otherwise we choose one of the
points
above
$Q_j$ arbitrarily as $Q_j^{+}$ and the other as $Q_j^{-}$. Then we can
associate to
$e_j$ the line bundle ${\cal O}_{N}(Q_j^{+})\otimes {\cal
O}_{N}(-Q_j^{-})$ whose total degree is $0$, but whose multidegree is
different from
$0$ if $e_j$ is not a loop. By linearity we obtain a map
$$
u:C_1(\Gamma, \ZZ ) \longrightarrow \prod\limits_{i\in I}\mbox{Pic}(N_i).
$$
If we restrict this map to $H_1(\Gamma, \ZZ )$ then the images are line bundles
of multidegree $0$, and using the identification $JN={^tJ}N$ we thus have
a homomorphism
$$
c: H_1(\Gamma, \ZZ ) \longrightarrow JN={^tJ}N
$$
which we call the {\em classifying map} of the Jacobian $JC$.
Note that $JC$ is given by the classifying homomorphism via the negative
of the pushout.
It is not difficult to
check that this map corresponds to the extension (1) (see \cite [Proposition
7]{Nam1}). (Note that our definition of the classifying map depends on the
orientation of the graph $\Gamma$ and on the choice of $Q^+_j$ and $Q^-_j$ if
$e_j$ is a loop. Both choices also appear in the identification of the
kernel of
(1) with $H^1(\Gamma,k^{\ast})$ and correspond to an automorphism $z\mapsto
z^{-1}$ in some factor of the torus
$H^1(\Gamma, k^{\ast}$).)

Let $\iota:C\rightarrow C$ be an involution with quotient
$C'=C/\langle\iota\rangle$. We assume that $\iota$ is not the identity on
any of
the components $C_i$ of $C$. The quotient $C'$ is again a nodal curve. We
denote the quotient map by $\pi:C\rightarrow C'$. If $C$ (and hence $C'$) is
smooth, the norm map is defined by
$$
\begin{array}{rcl}
\pi_*=\Nm:\ JC & \rightarrow & JC'\\
{\cal O}(\Sigma n_i P_i) & \mapsto & {\cal O}(\Sigma n_i \pi(P_i)).
\end{array}
$$
Recall that there is a norm map
$$
\pi_*=\Nm : JC\ \rightarrow\ JC'.
$$
also in the case where $C$ is a nodal curve
(for a general definition see e.g. \cite[Proposition
II.6.5]{EGA}).
One can show easily that the norm map is given by
$$
\Nm : {\cal L} \mapsto ({\cal L}\otimes \iota^{\ast}{\cal L})/ \langle
\iota\rangle.
$$
This coincides with the norm map described above when $C$ is smooth.
\begin{lemma}\label{lem11}
The norm map $\Nm:\ JC \rightarrow JC' $ is
surjective.
\end{lemma}

\begin{Proof}
If ${\cal M}\in JC'$ then ${\cal L}:= \pi^{\ast}{\cal M}\in
JC$
and $\Nm ({ \cal L})={\cal M}^2$. Surjectivity then follows since $JC'$ is
2-divisible.
\hfill
\end{Proof}

We first want to understand the structure of the group
scheme
$$
P=\ker (\Nm:JC\rightarrow JC')_0,
$$
i.e. the component of the kernel of the norm map containing the identity.
We shall then ask about good \lq\lq compactifications\rq\rq of $P$.
As in the classical case, $P$ can be alternatively defined as
\begin{displaymath}
  P= \ker \left( (1+\iota):JC\to JC \right)_0
  = \im  ((1-\iota):JC\to JC).
\end{displaymath}

\begin{definition}
We call $P$ the {\em (open) Prym variety} associated to the double cover
$\pi:C\rightarrow C'$.
\end{definition}
We shall see later that $P$ carries the structure of a semiabelian variety.

We have 3 possible {\em types}
of double points, namely:
\begin{enumerate}
\item[(1)]
Fixed points of $\iota$ where the 2 branches are not exchanged. Then the image
under $\pi$ is again a node. We shall call these nodes {\em (branchwise)
fixed}.
\item[(2)]
Fixed points of $\iota$ where the 2 branches are exchanged. Then the image
of this node is a smooth point of $C'$ and we shall refer to such nodes as
{\em swapping} nodes.
\item [(3)]
Nodes which are exchanged under $\iota$. These we shall call
{\em non-fixed} nodes.
\end{enumerate}

Recall
that the group of Cartier divisors on $C$ is
$$
\mbox{Div} \,C=\bigoplus\limits_{x\in C_{\operatorname{reg}}}\ZZ x +
\bigoplus\limits_{Q \operatorname{singular}}K_Q^{\ast}/{\cal
O}_Q^{\ast}.
$$
After choosing local parameters for the 2 branches of $C$ which intersect at a
point $Q$ we have an identification (see also \cite[p.158]{B})
$$
K_Q^{\ast}/{\cal O}_Q^{\ast}\cong k^{\ast}\times \ZZ \times \ZZ .
$$
Here $K$ denotes the ring of rational functions on the
normalization $N$ and the two integers are the
multiplicities of the divisors on the branches. For the three types of
fixed points we have:
\begin{enumerate}
\item[(1)]
In this case
$\pi(Q)$ is
again a node. The involution $\iota^{\ast}$ and the map ${\pi}_{\ast}$
are given
by
$$
\begin{array}{lcl}
\iota^{\ast}(c, m, n) &=& ((-1)^{m+n} c, m, n),\\
{\pi}_{\ast}(c, m, n) &=& ((-1)^{m+n}c^2, m, n)).
\end{array}
$$
\item[(2)]
Here
$\pi(Q)$ is a
smooth point and we have
$$
\begin{array}{lcl}
\iota^{\ast}(c, m, n) &=& (\frac 1 c, n, m)\\
{\pi}_{\ast}(c,m,n) &=& (m+n)\cdot \pi(Q).
\end{array}
$$
\item[(3)]
If the nodes $Q$ and $Q'$ are interchanged, then
$$
\iota^{\ast}: K_{Q'}^{\ast}/{\cal O}_{Q'}^{\ast}\cong
K_Q^{\ast}/{\cal O}_Q^{\ast}
$$
and after identifying
$K_Q^{\ast}/{\cal O}_Q^{\ast}$
with $K_{\pi(Q)}^{\ast}/{\cal O}_{\pi(Q)}^{\ast}$ the map $\pi_{\ast}$
induces the identity.
\end{enumerate}

We shall denote the edges, resp. vertices of the graph $\Gamma'=\Gamma(C')$ by
$f_{j'}, j'\in J'$ resp. $w_{i'}, i'\in I'$. The involution $\iota$ acts on
the sets $I$ and $J$ and hence by linearity on
$C_0(\Gamma,\ZZ)$ and
$C_1(\Gamma,\ZZ)$ and we can choose orientations of the graphs $\Gamma$ and
$\Gamma'$ which are compatible with $\iota$ in the following sense: If $e_j$ is
an edge in $\Gamma$ from $v_i$ to $v_{i'}$, $i\neq i'$ with $e_j\neq
e_{\iota(j)}$, then $e_{\iota(j)}$ goes from $v_{\iota(i)}$ to $v_{\iota(i')}$.
Otherwise we choose an arbitrary orientation on $e_j$. If $f_{j'}$ is an
edge of
$\Gamma '$ which is not a loop, it is either the image of
an edge
$e_j$ from $v_i$ to
$v_{i'}$ where $\iota$ fixes $j, i$ and $i'$ or
it is the image of two edges $e_j$
and $e_{\iota(j)}$ with $j\neq\iota(j)$. In the first case we orient
$f_{j'}$
in such a way that it goes from $w_{\pi(i)}$ to $w_{\pi(i')}$. In the
second case we choose the orientation in such a way that $f_{j'}=f_{\pi(j)}$
goes from $w_{\pi(i)}$ to $w_{\pi(i')}$ if $e_j$ goes from $v_i$ to
$v_{i'}$. By our choice of the orientation on $\Gamma$ this is well defined.
For a loop in $\Gamma '$ we choose an arbitrary orientation.

We define homomorphisms
$$
\iota_l:C_l(\Gamma,\ZZ)\rightarrow C_l(\Gamma,\ZZ),\ l=0,1.
$$
For $l=0$ this is defined by $\iota_0(v_i)=v_{\iota(i)}$.
For $l=1$ we set
$\iota_1(e_j)=\pm e_{\iota(j)}$ where we always use the $+$-sign with the
following two exceptions: we set $\iota_1(e_j)=-e_j$ if either $e_j$ comes
from a loop which corresponds to a swapping node or if
$e_{j}$ is an edge from $v_i$ to $v_{i'}, i\neq i'$ and
$j=\iota(j), i'=\iota(i)$.
It follows directly
from the definitions that $\iota_0$ and $\iota_1$ commute with the boundary
map $\partial$ and dually with the coboundary map $\delta$. Hence $\iota$
induces involutions
$$
\iota: H_1(\Gamma, \ZZ)\rightarrow H_1(\Gamma, \ZZ), \quad
\iota: H^1(\Gamma, \ZZ)\rightarrow H^1(\Gamma, \ZZ).
$$
The involution $\iota:C\rightarrow C$ also induces an involution on the
extension
$$
1\rightarrow T_C=H^1(\Gamma, \ZZ)\otimes k^{\ast}\rightarrow JC\rightarrow
JN\rightarrow 0.
$$
\begin{lemma}\label{lem12}
Under the identification $T_C=H^1(\Gamma, \ZZ)\otimes k^{\ast}$ the
restriction
of the involution $\iota:JC\rightarrow JC$ to $T_C$ is induced by
$\iota:H^1(\Gamma, \ZZ)\rightarrow H^1(\Gamma, \ZZ)$.
\end{lemma}
\begin{Proof}
Recall from the proof of \cite [Proposition 10.2]{OS} that we have an exact
sequence
$$
1\rightarrow {\cal O}_C^{\ast}\rightarrow {\cal
O}_N^{\ast}\stackrel{\alpha}{\rightarrow}\bigoplus\limits_{j\in J}
k_{Q_j}^{\ast}\rightarrow 1
$$
where $k ^{\ast}_{Q_j}$ is the skyscraper sheaf with fibre $k^{\ast}$
at the
node
$Q_j$. (The identification of the stalk of the cokernel over a node $Q_j$
with $k^{\ast}$ depends on choosing an order of the two branches
at $Q_j$. Whenever $Q_j$ defines an edge in $\Gamma$ which is not a loop then
we shall choose this order according to the chosen orientation of $\Gamma$.
If $Q_j$ gives rise to a loop then we can choose the order arbitrarily.)
The involution $\iota$ acts on ${\cal O}_C^{\ast}$ and on ${\cal
O}_N^{\ast}$ and hence also on the quotient. If $Q_j$
is a branchwise fixed node, then $\iota$ acts
trivially on $k_{Q_j}^{\ast}$, if $Q_j$
is a swapping node,
then $\iota$ acts by $z\mapsto z^{-1}$ and otherwise $\iota$
interchanges
$k_{Q_j}^{\ast}$ and $k^{\ast}_{Q_{\iota(j)}}$. The above sequence
induces a
long exact sequence
$$
1\rightarrow H^0({\cal O}_C^{\ast})\rightarrow H^0({\cal
O}_N^{\ast})\stackrel{\alpha}{\rightarrow}\bigoplus\limits_{j\in J}
k_{Q_j}^{\ast}
\rightarrow JC\rightarrow JN\rightarrow 0.
$$
The extension (1) is then a consequence of the observation that the map
$\alpha:H^0({\cal O}_N^{\ast})\rightarrow \oplus_{j\in
J}k^{\ast}_{Q_j}$
coincides with the coboundary map $\delta:C_0(\Gamma,k^{\ast})\rightarrow
C_1(\Gamma,k^{\ast})$. Under this identification the action of $\iota$
coincides with the action induced by $\iota_0$ and $\iota_1$
and this proves the claim.\hfill
\end{Proof}

The norm map $\Nm:JC\rightarrow JC'$ induces a diagram
$$
\diagram
1\rto &T_C=H^1(\Gamma,k^{\ast})\dto^{\Nm_T}\rto
&JC\dto^{\Nm}\rto & JN\dto^{\Nm}\rto & 0\\
1 \rto & T_{C'}=H^1(\Gamma',k^{\ast}) \rto & JC' \rto & JN'\rto & 0.
\enddiagram
\qquad (\mbox{D1})
$$
Note that all vertical maps are surjective. As before this follows since
$\Nm(\pi^{\ast}{\cal M})={\cal M}^2$ for ${\cal M}\in JC'$ and since the
groups involved are 2-divisible.

Adding the kernels of the vertical maps in diagram (D1) we obtain the following
diagram with exact rows and columns:

$$
\diagram
&1\dto &0\dto &0\dto&\\
1 \rto & T'_P \dto \rto &P' \dto \rto & K'_N \dto \rto & 0\\
1 \rto & T_C \dto \rto &JC\dto\rto &
JN\dto\rto & 0\\
1 \rto & T_{C'} \dto \rto &JC'\dto\rto
&JN'\dto\rto &0\\
&1 &0 &0&
\enddiagram
\qquad(\mbox{D}2)
$$
None of the group schemes in the top row need be connected. We had already
defined $P$ as the identity component of $P'$. Similarly we denote the
identity components of the other groups in the top row by
$$
K_N=(K'_N)_0,\quad T_P=(T'_P)_0.
$$
The group scheme $K_N$ is an abelian variety
which is the Prym variety for the double cover $N\to N'$,
and we have a morphism
$P\rightarrow K_N$ whose fibers, however, need not be connected. The identity
component of the kernel $K_P$ of this map is isomorphic to $T_P$ and we have a
commutative diagram
$$
\begin{array}{ccccccccc}
&&&&0&&0&&\\
&&&&\downarrow&&\downarrow&&\\
&&&&T_P&\hookrightarrow&K_P&&\\
&&&&\downarrow&&\downarrow&&\\
&&&&P&=&P&&\\
&&&&\downarrow&&\downarrow&&\\
0&\rightarrow&H&\rightarrow&A&\rightarrow&K_N&\rightarrow&0\\
&&&&\downarrow&&\downarrow&&\\
&&&&0&&0&&
\end{array}
$$
where $H\cong K_P/T_P$ is a finite group. In particular $A$ is an abelian
variety and the left hand column gives $P$ the structure of a semiabelian
variety. For future use we have to determine the classifying homomorphism
of this semiabelian variety. Before we do this we have to describe the
various lattices which will play a role.

The involution $\iota$ acting on the lattices $H^1(\Gamma, \ZZ)$, resp.
$H_1(\Gamma,\ZZ )$ defines
($\pm) 1$-eigenspaces $[H^1(\Gamma,\ZZ)]^{\pm}$, resp. $[H_1(\Gamma,
\ZZ)]^{\pm}$.

\begin{proposition}\label{prop14}
$T_P=[H^1(\Gamma, \ZZ)]^{-}\otimes k^{\ast}$.
\end{proposition}

\begin{Proof}
We first observe that $T_P\subset \mbox{Ker}
((1+\iota)\otimes\mbox{id}_{k^{\ast}})$. This follows since $\Nm({\cal
L})=({\cal L}\otimes\iota^{\ast}{\cal L})/{\iota}={\cal O}_{C'}$, implies that
${\cal L}\otimes\iota^{\ast}{\cal L}={\cal O}_C$. The identity component of
$\mbox{Ker}((1+\iota)\otimes\mbox{ id}_{k^{\ast}})$ is
$\mbox{Ker}(1+\iota)\otimes k^{\ast}=[H^1(\Gamma,\ZZ)]^{-}\otimes k^{\ast}$,
hence $T_P\subset [H^1(\Gamma,\ZZ)]^{-}\otimes k^{\ast}$. Next we claim
that for
every ${\cal L}$ with ${\cal L}\otimes\iota^{\ast}{\cal L}={\cal O}_C$ the
image
${\cal M}=\Nm({\cal L})$ is 2-torsion. This follows since ${\cal
M}^2=\Nm(\pi^{\ast}{\cal M})=\Nm({\cal L}\otimes\iota^{\ast}{\cal
L})=\Nm({\cal O}_C)={\cal O}_{C'}$. Since the 2-torsion elements in $JC'$ are
discrete the inclusion $[H^1(\Gamma, \ZZ)]^{-}\otimes k^{\ast}\subset T_P$
follows if there is one element in $[H^1(\Gamma,\ZZ)]^{-}\otimes k^{\ast}$
whose
image under the norm map is ${\cal O}_{C'}$. Since ${\cal O}_C$ is in
$[H^1(\Gamma,\ZZ)]^{-}\otimes
k^{\ast}$ we are done.\hfill
\end{Proof}
\begin{remark}
The above proof also shows that
$K_P=[H^1(\Gamma,\ZZ)]^{-}\otimes k^{\ast}\times (\mu_2)^N$ for some
integer $N$, where $\mu_2=\{\pm1\}$.
\end{remark}

For future use we fix the notation
$$
X=H_1(\Gamma,\ZZ),\ [X]^-=[H_1(\Gamma,\ZZ)]^-.
$$
Note that the dual lattice $X^{\ast}$ is the group of $1$-parameter
subgroups of the torus $T_C$, i.e.
$T_C=X^{\ast}\otimes  k^{\ast}$.

Recall (cf.\cite[Theorem (74.3)]{CR}) that every involution $\iota$ on a
lattice of finite rank
has, with respect to a suitable basis, the form
$$
\iota = k (1) + l (-1)+m\begin{pmatrix}
0 & 1\\ 1 & 0\end{pmatrix}.
$$
Such a decomposition is not canonical, but the integers $k,l$ and $m$ are
independent of the given decomposition. We consider the homomorphism
$$
\begin{array}{rcl}
\pi^{-}: H_1(\Gamma, \ZZ) & \longrightarrow & H_1(\Gamma,\frac 12 \ZZ)\\
h & \longmapsto & \frac 12 (h- \iota (h))
\end{array}
$$
and define
$$
X^-=\pi^{-}(H_1(\Gamma, \ZZ)).
$$
(Note that $X^-$ is not the $(-1)$-eigenspace of the lattice $X$.) Then
there is
an exact sequence
$$
0\rightarrow [H_1(\Gamma, \ZZ)]^{+}\rightarrow H_1(\Gamma,
\ZZ) \stackrel{{\pi}^-}{\rightarrow} X^- \rightarrow 0
$$
and the quotient map ${\pi}^-:H_1(\Gamma, \ZZ)\rightarrow X^-$ is
the dual to the
inclusion $[H^1(\Gamma, \ZZ)]^{-}\subset H^1(\Gamma, \ZZ)$.
Also note that after tensoring with $\QQ$ the map $\pi^-$ identifies
the $(-1)$-eigenspace $H_1(\Gamma, \RR)^-$
with $X^-_{\RR}=X^-\otimes \RR $. Over the rationals ${\pi}^-$ is nothing
but the projection onto the $(-1)$-eigenspace.
Via the map ${\pi}^-$ \
we can consider $[X]^-$ as a sublattice of $X^-$ with
$$
X^-/[X]^- \cong ( \ZZ/2)^m.
$$
Finally note that since $X^-$ is dual to $[H^1(\Gamma, \ZZ)]^-$, it is
the character group of the torus $T_P$ (cf Proposition \ref{prop14}).

The abelian base of the semiabelian variety $P$ is $A$ and the character
group of its toric part is $X^-$. By the general theory of semiabelian
varieties the group scheme $P$ is defined by its classifying homomorphism
$c^-:X^-\rightarrow {^tA}$. We have to determine this homomorphism.
In order to do this recall the following facts about homomorphisms of
semiabelian varieties.
Assume that $1\to T_i \to G_i \to A_i \to 0$;
$i=1,2$ are two semiabelian varieties with  classifying
homomorphisms $c_i:X_i\to {^tA}_i$. Then
giving a homomorphism $\varphi:G_1\to G_2$ is equivalent to giving two
homomorphisms ${^t\varphi}_T,{^t\varphi}_A$ making the following diagram
commutative:
$$
\diagram
X_2\dto^{c_2}\rto^{{^t\varphi}_T} &X_1\dto^{c_1}\\
{^tA}_2 \rto^{{^t\varphi}_A} & {^tA}_1.
\enddiagram
$$

The following is a standard lemma whose proof we omit. In this lemma
we assume for simplicity that the characteristic of the field $k$
does not divide the order of any of the finite groups appearing. Note that
we will be in this position since we have assumed that the characteristic
of $k$ is different from $2$ and since the finite groups which will appear
are all $2$-groups.

\begin{lemma}\label{lemsemiabelian}
Let $\varphi:G_1\rightarrow G_2$ be a homomorphism of semiabelian varieties.
  \begin{enumerate}
  \item
The identity component of the kernel
$(\ker\varphi)_0$ is a semiabelian variety defined by the last
    column in the diagram
    \begin{displaymath}
      \diagram
      X_2\dto^{c_2}\rto^{{^t\varphi}_T}
      & X_1\dto^{c_1}\rto
      & \coker{{^t\varphi}_T}\dto^{c_1'}\rto
      & \coker{{^t\varphi}_T/{\rm Torsion}}\dto
      \\
      {^tA}_2 \rto^{{^t\varphi}_A}
      & {^tA}_1\rto
      & \coker{{^t\varphi}_A}\rto
      & \coker{{^t\varphi}_A}/c'_1\left(
      {\rm Torsion}(\coker{{^t\varphi}_T})\right).
      \enddiagram
    \end{displaymath}
 \item The image $\im\varphi$ is a semiabelian variety defined by the
third column in the diagram
\end{enumerate}
\vspace{-3mm}
{\small $$
\diagram
      \ker{^t\varphi}_T \rto\dto^{c_2}
      & X_2\rto^{{^t\varphi}_T}\dto^{c_2}
      & \im{^t\varphi}_T\dto\rto
      & X_1\dto^{c_1}\\
      ((\ker{^t\varphi}_A)_0+c_2(\ker{^t\varphi}_T))\rto
      & {^tA}_2\rto
      & {^tA}_2/((\ker{^t\varphi}_A)_0+c_2(\ker{^t\varphi}_T))\rto
      & {^tA}_1.
\enddiagram
$$}

    Here $(\ker{^t\varphi}_A)_0\subset ((\ker{^t\varphi}_A)_0+
    c_2(\ker{^t\varphi}_T))
    \subset \ker{^t\varphi}_A$ with finite co\-kernels.
 \end{lemma}

We can now obtain a description of the classifying
homomorphism of $P$.
For this we apply the first part of this lemma to
$P=\left(\ker(1+\iota)\right)_0$. Then we have
\begin{eqnarray}
  \label{eq:lattice_for_Prym}
X\overset{1+\iota}{\longrightarrow}
X\to X/(1+\iota)(X) \to (X/(1+\iota)(X))/{\rm Torsion}=X^-
\end{eqnarray}
and the torsion group is identified with the image of $X^+$ and,
moreover, is isomorphic to $(\ZZ/2)^k$, where $k$ is the number of
$(+1)$-blocks in the decomposition of the involution~$\iota$. This
immediately gives

\begin{proposition}\label{propclasshom}
The classifying homomorphism of the semiabelian variety $P$ is given by
the right hand vertical map of the following diagram
$$
\diagram
X\dto^{c}\rto &X^-\dto^{c^-}\\
\tu{JN} \rto & \tu{A} = \tu{K_N}/ \im(\ZZ/2)^k.
\enddiagram
$$
The abelian variety $A$ is a $2^a$ cover of the Prym variety $K_N$ of the
normalization, where $a\leq k$ and $k$ is the number of $(+1)$-blocks in
the decomposition of the involution $\iota : X \to X$.
\end{proposition}

\subsection{Geometric properties of the kernel of the norm map}
in this subsection we shall study the geometric properties of the kernel of
the norm map. In particular, we shall determine the number of its connected
components
and the type of the induced polarization on the abelian part $A$ of $P$.
We write the set $S$ of nodes of $C$ as
$$
S=S_e\cup\iota(S_e)\cup S_f \cup S'_f
$$
where $S_e \cup \iota (S_e)$ is the set of nodes which are exchanged under
$\iota$, $S_f$ is the set of nodes which are branchwise fixed
and $S'_f$ is the set of swapping nodes. For every
node $s$ we have (after a choice of local parameters) a homomorphism
$$\begin{array}{rcl}
k_s^{\ast}&\longrightarrow& JC\\
c_s&\longmapsto&{\cal O}(c_s, 0, 0).
\end{array}
$$
(Note that we are not claiming that this map is injective. If $s$ is a
disconnecting node, then the above homomorphism is the constant homomorphism
mapping $k_s^{\ast}$ to the trivial line bundle.) Recall that $P'$ is the
kernel of the norm map.

\begin{proposition}\label{prop15}
Every element ${\cal L}\in P'$ is of the form
$$
{\cal L}\cong\bigotimes\limits_{s\in S'_f}{\cal O}(c_s, 0, 0)\otimes {\cal
M}\otimes \iota^{\ast}{\cal M}^{-1}
$$
where ${\cal M}\in
%VA010800 J
\Pic
%ENDVA010800
(C)$ is a line bundle on $C$
whose degree on every
component of $C$ is $0$ or $1$.
\end{proposition}

\begin{Proof}
As in \cite[Lemma 1]{M} or \cite[Lemma (3.3)]{B} we can assume that ${\cal
L}={\cal O}(D)$ where $D$ is a Cartier divisor with ${\pi}_{\ast}(D)=0$.
For branchwise fixed nodes  we have
$$
(-1, 0, 0) = (1, 0, 1)-\iota^{\ast}(1, 0, 1)
$$
and for swapping nodes  we have
$$
(1, 1, -1) = (1, 1, 0) -\iota^{\ast}(1, 1, 0).
$$
Hence we have a decomposition
$$
D = D'-\iota^{\ast} D'+\Sigma_{s\in S'_f}(c_s, 0, 0).
$$
This shows the claim apart from the assertion on the degrees of ${\cal M}$.
We can
always replace ${\cal M}$ by ${\cal M}\otimes{\pi}^{\ast}{\cal N}$
for some line
bundle ${\cal N}$ on $C'$. This shows immediately that we can assume that the
degree of ${\cal M}$ is $0$ or $1$ on components $C_i$ of $C$ which are fixed
under
$\iota$. Now assume that $C_i$ and $C_{\iota(i)}$ are two components which are
interchanged. Since ${\cal L}\in P'\subset JC$
we have $\deg {\cal L}|_{C_i}=\deg {\cal L}|_{C_{\iota(i)}}=0$
and this implies that
$\deg
{\cal M}|_{C_i}=\deg {\cal M}|_{C_{\iota(i)}}$. But then, after replacing
${\cal M}$ by ${\cal M}\otimes \pi^{\ast}{\cal N}$ for some line bundle
${\cal N}$ on $C'$
we can in fact assume that $\deg {\cal M}|_{C_i}=0$ for all
components $C_i$
which are not fixed under the involution $\iota$.
\hfill \end{Proof}

In order to compute the number of components of $P'$ we have to determine
in how
far we can normalize the multidegree of the line
bundles ${\cal M}$ in Proposition \ref{prop15}. The proof of this proposition
shows that we can assume that the degree of ${\cal M}$ is $0$ on components
of $C$
which are not fixed under $\iota$. As before we denote by $B$ the union of the
components of $C$ which are fixed by $\iota$. Let $N'_B$ be the partial
normalization of $B$ obtained by blowing up the nodes of $B$ which are not
branchwise fixed. Then $N'_B$ decomposes into connected components
$$
N'_B=N^1_B\cup\ldots\cup N_B^{n_B}\cup L^1_B \cup\ldots\cup L_B^{l_B}
$$
where the components $N_B^i$ have no smooth fixed points with respect to
$\iota$, whereas the components $L_B^j$ have such a fixed point.

\begin{proposition}\label{prop16}
The number of components of the kernel $P'$ of the norm map
equals $2^{n_B}$ where $n_B$ is the
number of connected components of the curve $N'_B$ on which the involution
$\iota$ has no smooth fixed points.
\end{proposition}
\begin{Proof}
We first show that the number of components of $P'$ is at most
$2^{n_B}$. The involution $\iota$ acts on each of the components
$N_B^i$, $L_B^j$ defining varieties $P_{N^i_B}$, $Q_{L^j_B}$
which are defined as the kernel of the restricted norm maps.
It was shown by Beauville \cite[Lemma 3.3]{B} that $P_{N^i_B}$
has $2$ components. In particular Beauville showed that one can
normalize the multidegree of ${\cal M}|_{N_B^i}$ to $(0,0,\dots,0)$ or
$(1,0,\dots,0)$ depending on whether the sum of the degrees is even or odd.
The curves $L_B^j$ have a smooth fixed point $x=\iota(x)$.
Since ${\cal O}(x)\otimes{\iota}^{\ast}{({\cal O}(x))}^{-1}$ is trivial
we can change the parity of the multidegree of ${\cal M}$ in this case
and hence in this case $Q_{L^j_B}$ is irreducible. We had already remarked that
we can assume that the multidegree on $D\cup\iota(D)$ is $(0,0,\dots,0)$.
This shows that it suffices to consider at most $2^{n_B}$ possible
multidegrees.
Since the variety of line bundles on $C$ with fixed multidegree is
irreducible, it
follows from Proposition \ref{prop15} that the number of components
of $P'$ is at most
$2^{n_B}$.

To prove that the number of components is at least $2^{n_B}$ it is
sufficient to show that
the restriction $P'\rightarrow \prod\limits_i P_{N^i_B}\times
\prod\limits_j Q_{L^j_B}$ is surjective. We choose line
bundles ${\cal M}_i\in P_{N^i_B}$, resp. ${\cal M}_j\in Q_{L^j_B}$ and the
trivial
bundle ${\cal O}_{D\cup\iota(D)}$. We want to glue these line bundles to a line
bundle ${\cal M}$ on $C$ which is contained in $P'$. To do this we have to
specify the gluing at the nodes which are not branchwise fixed. The gluing over
the swapping nodes can be chosen arbitrarily. Since we can
always choose the gluing at a node $s$ and its image $\iota(s)$ in such a way
that $\Nm({\cal M})={\cal O}_C$ the claim follows.
\hfill
\end{Proof}

Our next aim is to determine the type of the polarization on $A$ which is
induced from the principal polarization on $JN$ via the
map $A\rightarrow K_N\rightarrow JN$.
Roughly speaking $A$ is made up from $3$ types of building blocks which
arise as follows:
\begin{enumerate}
\item [(1)]
Let $D_i$ and $D_{\iota(i)}$ be $2$ components of $C$ which are exchanged
under the involution $\iota$ and denote the normalizations of these
components by $N_{D_i}$ and $N_{D_{\iota (i)}}$. Then
the kernel of the norm map defined by the quotient map $N_{D_i}\cup
N_{D_{\iota(i)}}\rightarrow N_{D_i}=N_{D_{\iota(i)}}$
is $JD_{N_i}$ embedded as the
anti-diagonal in
$JN_{D_i}\times JN_{D_i}$. The restriction of the product polarization on
$JN_{D_i}\times JN_{D_i}$ to the (anti-)diagonal is twice a
principal polarization.
\item [(2)]
The involution $\iota$ on a component $N_B^i$ defines an abelian variety
$P_{N_B^i}$ and the induced polarization is twice a principal polarization.
This was shown by Beauville \cite{B}.
\item [(3)]
The third building block comes from the components $L_B^j$. To simplify the
notation we consider a connected nodal curve $L$ together with an involution
$\iota$ where we assume that all nodes are branchwise fixed.
We do, however, allow that
$\iota$ has smooth fixed points. Let $L'=L/\langle \iota\rangle$ and we denote
by $N_L$ the normalization of $L$, resp. by $N'_{L'}$ that of $L'$. Then we
have varieties
$$
\begin{array}{lclcl}
Q & = & \mbox{Ker} (\Nm:JL &\rightarrow & JL'),\\
R & = & \mbox{Ker} (\Nm:JN_L & \rightarrow & JN'_{L'}).
\end{array}
$$
We have already observed that $Q$ and $R$ are connected. Since all nodes are
branchwise fixed the map
$$
Q\rightarrow R\subset JN_L
$$
is finite and hence $Q$ is an abelian variety.
\end{enumerate}

\begin{lemma}\label{lem17}
Assume that $L$ is a connected nodal curve with an
involution $\iota:L\rightarrow L$ which has only branchwise fixed nodes
and let $r$ be the number of smooth fixed points of $\iota$.
Let $Q$ be the abelian variety defined above. Then
\begin{enumerate}
\item[{\rm(i)}]
The induced polarization on $Q$ is of type $(1,\ldots,1,\ 2,\ldots,2)$.
\item[{\rm(ii)}]
It is of type $(2,\ldots,2)$ if and only if $r=0$ or $2$.
\end{enumerate}
\end{lemma}

\begin{Proof}
The number of smooth fixed points of $\iota$ is even. By
identifying these points pairwise we obtain a curve $\overline{C}$ together
with an involution $\overline{\iota}$ which is of Beauville type, i.e. all
double points are branchwise fixed and there are no smooth fixed points. We can
consider $C$ as a partial normalization of $\overline{C}$. The claim then
follows from \cite[Lemma 1]{DL}. \hfill
\end{Proof}

\begin{remark}
The number of entries equal to $1$ in the above lemma is $r'-1$ if there are
$r=2r' > 0$ smooth fixed points.
\end{remark}

We can now prove the
\begin{proposition}\label{prop18}
The polarization of the abelian variety $A$ which is induced from the map
$A\rightarrow JN$ is always of type $(1,\ldots,1,2,\ldots,2)$. It
is twice a principal polarization if and only if the
following condition ($\dagger$) holds:\\
($\dagger$) Every connected component of the partial
blow-up $N'_B$ has at most $2$ smooth fixed points with respect to $\iota$.
\end{proposition}
\begin{Proof}
We consider
$$
A':=JN_{D_1}\times\ldots \times JN_{D_{n_e}}
\times P^0_{N^1_B}\times\ldots\times
P^0_{N^{n_B}_B}\times Q_{L^1_B}\times\ldots\times Q_{L^{l_B}_B}
$$
where $P^0_{N^k_B}$ is the identity component of $P_{N^k_B}$. By
Lemma (\ref{lem17}) and the discussion preceding this lemma the polarization
induced
by the map $A'\rightarrow JN$ is twice a principal polarization if and only
if ($\dagger$) holds. We claim that $A'=A$. Pulling back line bundles to the
partial normalization defined by blowing up all nodes which are not
branchwise fixed,
resp. to
the normalization $N$ defines a diagram:
$$
\diagram
P \rrto \drto & & K_N. \\
& A'\urto
\enddiagram
$$
The map $A'\rightarrow K_N$ is finite.
In order to show that $A'=A$ it is,
therefore, enough to prove that the kernel of
the map $P\rightarrow A'$ is connected. Let
${\cal L}$ be an element in this kernel. Then we can
write ${\cal L}={\cal O}_C(D)$
with $D=D_1+D_2+D_3$ where the support of the divisor $D_k$ is contained in
the set of nodes of type
$(k)$. Again we can argue as in \cite[Lemma 1]{M} to conclude that we may
assume that $\pi_{\ast}(D)=0$. Moreover we claim that we may assume $D_1=0$.
Indeed, the restriction of ${\cal O}(D)$ to each of the connected components
$N^i_B$ or $L^j_B$ is trivial. Hence there exists for each of these components
a rational function $g^i$ or $h^j$ such that the restriction of $D$ to $N^i_B$
or $L^j_B$ is equal to the principal divisor $(g^i)$ or $(h^j)$. We can assume
that $g^i$ and $h^j$ are $\pm 1$ on each irreducible component of $N^i_B$
and
$L^j_B$. We extend these rational functions to a rational function on $C$ by
setting it equal to $1$ for every component which is not fixed under $\iota$.
The divisor
$$
D':=D-\sum_i(g^i)-\sum_j (h^j)\sim D
$$
is only supported on nodes which are either swapping or non-fixed nodes.
We claim that we have still
$\pi_{\ast}(D')=0$. There is nothing to check for  swapping nodes.
Exchanged nodes come in pairs, say $\{Q,Q'\}$, which
are interchanged by $\iota$ and there are the following possibilities.
The first possibility is that $\{Q,Q'\}$ is contained in $1$ or $2$ curves
which are invariant under $\iota$. The second
possibility is that $\{Q,Q'\}$ is contained in the
intersection of a curve $B_i$ with a curve $D_j$ and its image $D_{\iota(j)}$
under the involution $\iota$.
The third possibility is that none of the curves containing $\{Q,Q'\}$
is fixed under $\iota$.
In either case the property $\pi_{\ast}(D)=0$ is
not effected by changing $D$ to $D'$. We claim that the line bundles ${\cal
O}(D')$ with $\pi_{\ast}(D')=0$ and $D'$ supported at swapping nodes
and exchanged nodes such
that the pull back of $D'$ to $N$ is trivial are parameterized by an
irreducible
variety. This follows since the variety of line bundles of the form
$$
\bigotimes_{s\in S'_f}{\cal O}(c_s,0,0)\otimes\bigotimes_{s\in S_e}{\cal
O}(d_s,0,0)\otimes
\bigotimes_{s\in S_e}{\cal O}(d^{-1}_{\iota(s)}, 0,0)
$$
is the image of a finite product of tori $k^{\ast}$ and hence
irreducible.\hfill
\end{Proof}

\begin{remark}
The description of $A$ in the proof of the above theorem together with
the remark after
Lemma \ref{lem17} allows us to compute the degree of the induced
polarization on $A$.
\end{remark}

\begin{remark}
  Geometrically the most interesting case occurs when $(C, \iota)$ can
  be smoothed to a curve $\left(C(t), \iota(t)\right)$, where $\iota(t)$
  has $0$ or $2$ smooth fixed points. This is the case if $(C, \iota)$
  has either no smooth fixed points and at most one swapping node
  or $2$ smooth fixed points and no swapping nodes. In either case
  the above proposition shows that the induced polarization is twice a
  principal polarization, i.e. of type $(2,\ldots, 2)$.
\end{remark}

%%%%%%%%%%%%%%%%%%%%%%%%%%%%%%%%%%%

\section{Theory of degenerations}
\label{sec:Theory of degenerations}
\subsection{General theory of degenerations of abelian varieties}
\label{sub:General theory of degenerations of abelian varieties}

For the reader's convenience we shall review the basic facts about
degenerations of abelian varieties in a form relevant to this article.
Unfortunately, a certain amount of technicalities is unavoidable in
this subject.
We shall begin with the Mumford-Faltings-Chai uniformization of
abelian varieties over the quotient field of a complete normal ring,
\cite[Ch.II]{FC}.

The general setup is that $R$ is a noetherian normal integral domain
complete w.r.t. an ideal $I=\sqrt{I}$, which is a completion of a normal
excellent ring. Let $S=\Spec R$, $S_0=\Spec R/I$ and let
$\eta$ be the generic point of
$S$ and $K$ be the fraction field of $R$.
Assume we are given the following data:

\begin{enumerate}
\item[(a)] Let $G/S$ be a semiabelian scheme whose generic fiber $G_{\eta}$
  is abelian such that $G_{S_0}$ is an extension of an abelian scheme
  $A_0$ by a split torus~$T_0$,
\item[(b)] Let $\calL$ be an invertible sheaf on $G$ rigidified at the zero
  section such that $\calL_{\eta}$ is ample.
\end{enumerate}

The main result of \cite[Ch.II]{FC} is that the pair $(G,\calL)$ is
equivalent to a set of \emph{degeneration data} (d0) - (d4)
which we will describe
below. Here small letters will indicate that we describe the degeneration
data of a semiabelian scheme over the base $S$ and capital letters
will later indicate that we describe the degeneration data of the central
fiber over $S_0$. The degeneration data of the pair $(G,\cal L)$
given above are as follows:

\begin{enumerate}
\item[(d0)] An abelian scheme $A/S$ with an ample rigidified sheaf
  $\calM$. This sheaf determines a polarization $\lambda:A\to \tu A$,
  where $\tu A$ is the dual abelian scheme over $S$.
\item[(d1)]
  \begin{enumerate}
  \item A semiabelian scheme, otherwise known as Raynaud extension
    $$1\to T\to \wG \to A\to 0$$
    over $S$ with a torus $T$. Since the torus $T_0$ is split, the
    character group $X$ of $T$ is a
    constant group scheme over $S$, i.e. $X \cong \ZZ^r$.
    Then $\wG$
    corresponds uniquely to a homomorphism $c:X\to \tu A$ via the
    negative of the pushout.
  \item A second semiabelian scheme $$1\to \tu T\to \tu \wG \to \tu A\to
    0$$ with a split torus part, given by a homomorphism $\tu c:Y\to
    A$.
  \end{enumerate}
\item[(d2)] An inclusion of lattices $\phi:Y\to X$ with a finite
cokernel such that $c\circ\phi=\lambda\circ\tu c$. This corresponds to
a homomorphism $\wG\to \tu \wG$.
\item[(d3)] A bihomomorphism $\tau:Y\times X\to (\tu c\times
  c)^*\calP\inv_{\eta}$, where $\calP_{\eta}$ is a Poincar\'e
  biextension on $(A\times \tu A)_{\eta}$, the universal Poincar\'e bundle
  with the zero section removed.
  Explicitly, this means that for every $y\in Y$,
  $x\in X$ we are given a nonzero section $\tau(y,x)$ of the sheaf
  $c(x)_{\eta}\inv$ over $\tu c(y)_{\eta}\simeq \Spec K$, and that
  $$\tau(y_1+y_2,x)\tau(y_1,x)\inv \tau(y_2,x)\inv =1$$
  and that a similar identity holds for the linearity in the second variable.
  The fact that $\calP$ is a biextension means that the left-hand side of
  the above identity is an element of the trivial $\GG_m$-torsor over
  $K$, so the identity makes sense.

  Here $\tau$ is required to satisfy the following positivity condition.
  The element $\tau(y_1,\phi(y_2))$ can be understood as a section of
  $(\id_A,\lambda)^*\calP_{\eta}\inv = \Lambda(\calM_{\eta}^*{}\inv)$ on the
  pullback by $(\tu c\times \tu c):Y\times Y \to (A\times A)_{\eta}$.
  Here $\calM^*$ is the $\GG_m$-torsor associated to the invertible
  sheaf $\calM$ and $\Lambda(\calM^*)=m^*\calM^* \otimes
  p_1^*\calM^*{}\inv \otimes p_2^*\calM^*{}\inv$ is the canonical
  symmetric biextension on $A_{\eta}{\times}A_{\eta}$, $p_i$ are the
  projections, and $m$ is the multiplication map. The positivity
  condition is that for all $y$ the element $\tau(y,\phi(y))$ is
  defined over the whole of $S=\Spec R$ and for every $n\ge0$ it is 0 modulo
  $I^n$ for all but finitely many $y$.

\item[(d4)] A cubical morphism $\psi:Y\to (\tu
  c)^*\calM_{\eta}^*{}\inv$. Explicitly,
  this means that for every $y\in Y$ we a have a nonzero section
  $\psi(y)$ of the sheaf $\calM_{\eta}\inv$ over $\tu c(y)_{\eta}$,
  and that
  $$\psi(y_1+y_2+y_3)
  \psi(y_2+y_3)\inv\psi(y_1+y_3)\inv\psi(y_1+y_2)\inv
  \psi(y_1)\psi(y_2)\psi(y_3) =1 $$
  (in particular, $\psi(0)=1$).
  Again, the fact that $\calM$ has a canonical cubical structure
  implies that the value of the left-hand-side is in the trivial
  $\GG_m$-torsor, and so the identity makes sense.

  Moreover, $\psi$ is required to satisfy the identity
  $$\psi(y_1+y_2)
  \psi(y_1)\inv \psi(y_2)\inv = \tau(y_1,\phi(y_2)) = \tau(y_2,
  \phi(y_1))$$
  on $(\tu c\times \tu c)^* \Lambda(\calM_{\eta}^*{}\inv)$. It follows that for
  all but finitely many $y$ the element $\psi(y)$ is a section of
  $\calM\inv$ over the whole of $S$, and, moreover, for each $n\ge 0$ all but
  finitely many $\psi(y)$ are 0 modulo $I^n$.
\end{enumerate}

The meaning of these degeneration data is as follows: $\tau(y,x)$ is
canonically an element of
$\Iso\left(T^*_{^t c(y)} c(x), c(x) \right)_{\eta}$, i.e. an invertible
element of
\begin{eqnarray*}
  & \Hom\left(T^*_{^t c(y)} c(x), c(x) \right)_{\eta} =
  \left( T^*_{^t c(y)} c(x)_{\eta} \right)^{-1}|_{\{0\}}
  \otimes c(x)_{\eta}|_{\{0\}} = \\
  & = \left( T^*_{^t c(y)} c(x)_{\eta} \right)^{-1}|_{\{0\}} =
  {\cal P}_{\eta}^{-1} \left( ^t c(y),c(x) \right)
\end{eqnarray*}
  (since $c(x)$ is linearized at 0),
and similarly $\psi(y)$ is an element of \linebreak
$\Iso\left(T^*_{c^t(y)}\calM,
\calM\otimes c(\phi(y)) \right)_{\eta},$
see \cite[p.44]{FC} for more details.
Therefore, $\tau$ describes an embedding $Y\to \wG_{\eta} = \Spec
\oplus_{x\in X} c(x)$ which via $\pi:\wG\to A$ sits over $\tu c:Y\to
A_{\eta}$, and $\psi$ describes an action of $Y$ on the sheaf
$\pi_{\eta}^*\calM_{\eta}$. The uniformization means that the abelian
variety $G_{\eta}$ can be thought of as the quotient $\wG_{\eta}/Y$
and the sheaf $\calL_{\eta}$  as the quotient
$(\pi_{\eta}^*{\calM}_{\eta})/Y$.

Moreover, by \cite[Thm.5.1]{FC} any section $s\in
\Gamma(G_{\eta},\calL_{\eta})$ can be written as a partial Fourier
series $s=\sum_{x\in X} \sigma_x(s)$, with $\sigma_x(s)\in
\Gamma\left(A,\calM\otimes c(x)\right)_{\eta}$ which satisfy the
following identity:
$$ \sigma_{x+\phi(y)}(s) = \psi(y) \cdot \tau(y,x) \cdot
T_{\tu c(y)}^* \sigma_x(s).$$
This makes sense since $\psi(y)\tau(y,x)$ is canonically an element of the
group
$\Iso\left(T^*_{c^t(y)} \left(\calM\otimes c(x) \right),
\calM\otimes c(x+\phi(y)) \right)_{\eta}.$

Note also that the degree of the polarization $\lambda_{\eta}:G_{\eta} \to
\tu G_{\eta}$ is the product of the degree of the polarization
$\lambda_A:A\to \tu A$ and $|X/Y|$. Moreover, one can show that there
is the following exact sequence:
$$ 0 \to \ker(\lambda_A)_{\eta} \to \ker\lambda_{\eta} / (Y^*/X^*)
\to (X/Y) \to 0,$$
where we denote $X^*= \operatorname{Hom}(X,\ZZ)$ and
$Y^*= \operatorname{Hom}(Y,\ZZ)$.

The degeneration data (d0)-(d4) give
corresponding degeneration data for the ``central fiber'' over
$S_0=\Spec R/I$. We shall later on describe how to reconstruct the
``central fiber'' from these data. To begin with the data (d0) - (d2) give
corresponding data (D0) - (D2) for an arbitrary base $S$. Namely we get

\begin{enumerate}
\item[(D0)] An abelian scheme $A_0/S_0$ with an ample rigidified sheaf
  $\calM_0$. This sheaf determines a polarization $\lambda_0:A_0\to \tu
  A_0$.
\item[(D1)]
  \begin{enumerate}
  \item A semiabelian scheme $$1\to T_0\to \wG_0 \to A_0\to 0$$
    with a
    split torus part, given by a homomorphism $c_0:X\to \tu A_0$.
  \item A second semiabelian scheme $$1\to \tu T_0\to \tu \wG_0 \to
    \tu A_0\to 0$$
    with a split torus part, given by a homomorphism
    $\tu c_0:Y\to A_0$.
  \end{enumerate}
\item[(D2)] An inclusion of lattices $\phi:Y\to X$ with a finite
cokernel such that $c_0\circ\phi=\lambda_0\circ\tu c_0$.
\end{enumerate}

We shall now restrict ourselves to the special situation of a
one-parameter degeneration. In this case $R$ is a complete discrete
valuation ring, $R/I=k$ is a field (not necessarily algebraically
closed) and $S_0=\Spec k$ is a one-point
scheme. In addition to the pair $(G_{\eta},\calL_{\eta})$ we shall also pick a
section $\theta_{\eta}\in\Gamma(G_{\eta},\calL_{\eta})$. In other words, we
shall look
at a pair $(G_{\eta},\Theta_{\eta})$, where $\Theta_{\eta}$ is the
Cartier divisor determined by $\theta_{\eta}$. Then by \cite[Thm.5.7.1]{A2}
this pair can, possibly after a base change, be extended uniquely to
a \emph{stable semiabelic pair} $(\wG\acts
\overline{P}\supset\Theta)$ over $S$.
This means in particular, that $\overline{P}_0$ is a
projective seminormal variety over $k$ with a $G_0$-action and
finitely many orbits, $\Theta_0$ is an ample Cartier
divisor on it not containing any orbits entirely,
and they satisfy certain conditions listed in\cite[1.1]{A2}. (In
\cite{A2} one starts with a $G_{\eta}$-torsor $P_{\eta}$, but we will
take $P_{\eta}=G_{\eta}$ here for simplicity). The central fiber
$(\overline{P}_0,\Theta_0)$ is described by the data (D0)-(D2)
and the additional
data (D3)-(D6) which we list below, following \cite{A2}.

Since $R$ is now a DVR and the ideal $I=(t)$ is principal, we can
write the bihomomorphism $\tau:Y\times X\to (\tu c\times
c)^*\calP\inv_{\eta}$ as
$$\tau(y,x) = t^{B(y,x)} \tau'(y,x),$$
where $B(y,x)$ is the valuation
of $\tau(y,x)$. This defines a bilinear form $B:Y\times X\to \ZZ$
such that $B|_{Y\times Y}$ is symmetric and positive definite. Considering
$\tau$ modulo the ideal $I$ we obtain the datum
\begin{enumerate}
\item[(D3)] A bihomomorphism
  $\tau_0=\tau\!\mod I:Y\times X\to (\tu c\times c)^*\calP\inv_{0}$
  on the central fiber.
\end{enumerate}

Applying the same procedure to the cubical morphism $\psi$ we get a
qua\-dra\-tic nonhomogeneous function $A:Y\to \ZZ$ such that $A(0)=0$,
whose quadratic part is $\frac{1}{2}B|_{Y\times Y}$. Taking $\psi$
modulo $I$ we obtain
\begin{enumerate}
\item[(D4)] A cubical morphism $\psi_0:Y\to (\tu
  c)^*\calM^*_{0}{}\inv$ on the central fiber.
\end{enumerate}

Finally, applying this to $\theta_{\eta}$ we obtain
\begin{enumerate}
\item[(D5)] For each $x\in X$, a section $\theta_{x,0}\in \Gamma(A_0,
  \calM\otimes c(x))$, such that
  $$\theta_{x+\phi(y),0} = \psi_0(y) \cdot \tau_0(y,x) \cdot
  T^*_{\tu c_0(y)}\theta_{x,0}. $$
\end{enumerate}

We also get a function $H:X\to\ZZ$,
$H(x)=\operatorname{val}_t(\theta_x)$. One easily checks the
following: The function $A:Y\to \ZZ$ extends uniquely to a quadratic
nonhomogeneous function $A:X\to\QQ$, and $H(x)=A(x) + r(x \!\mod
Y)$, where the last function depends only on the residue of $x$ modulo
$Y$. The function $H$ determines the last part of the data:
\begin{enumerate}
\item[(D6)] A cell decomposition $\Delta$ of the vector space
  $X\otimes\RR$ into polytopes. The vertices of these polytopes are in
  $X$, and it is periodic with respect to the lattice $Y$. The
  construction is as follows. $H$ is the ``height function'', i.e. one
  considers the convex hull of the countably many points $(x,H(x))$, $x\in
  X$. Since the leading term of $H$ is a positive definite quadratic
  form (whose first derivative is linear), it is easy to see that the
  lower envelope of this hull consists of countably many polytopes and
  that the projection of these polytopes to $X\otimes\RR$ gives a
  $Y$-periodic decomposition. This is the decomposition  $\Delta$.
\end{enumerate}

\begin{remark}
\begin{enumerate}
\item Only the equivalence classes of $\tau_0,\psi_0,\theta_0$ in the
  cohomology groups $H^1(\Delta,{\TT})$ and $Z^1(\widehat{\MM^*})$
  matter, see \cite{A2} for more details.
\item Changing the sheaf $\calM_{\eta}$ by a translated sheaf (hence,
  in the same polarization class), and the divisor $\Theta_{\eta}$ by
  the corresponding translated divisor, gives an isomorphic pair
  $(\overline{P}_0,\Theta_{0})$ if we do not choose
  an ``origin'' in $\overline{P}_0$.
\item The function $\tau_0$ (resp. $B$) plays the role of a 1-cohomology
  class, and the function $\psi_0$ (resp. $A$) the role of a
  1-cocycle restricting to this 1-cohomology class. So, changing
  $\psi_0$ without changing $\tau_0$ leads to an isomorphic pair
  $(\wG_0 \acts \overline{P}_0,\Theta_0)$,
  again, if we do not fix the ``origin''
  in~$\overline{P}_0$.
\item If the polarization $\lambda_{\eta}$ is principal, then the divisor
  $\Theta_{\eta}$ is superfluous, since for any sheaf $\calL_{\eta}$
  defining $\lambda_{\eta}$ the section $\theta_{\eta}$ is unique, up
  to a multiplicative constant.
\item In the case of principal polarization we must have $Y=X$, so for
  the decomposition $\Delta$ the two lattices of vertices and of
  periods coincide, and the periodic remainder function $r(x \!\mod
  Y)$ disappears. In this case, $\Delta$ is a classical Delaunay
  decomposition for the quadratic form $B$. By analogy, we will
  christen the decompositions appearing in the general case as
  \emph{semi-Delaunay decompositions}.
\item The whole construction can be repeated in entirely the same
  manner in the complex-analytic setting. The base $S$ then has to be
  replaced by the germ of a normal analytic space, and the ring $R$
  by the corresponding ring of convergent power series.
\end{enumerate}
\end{remark}

Vice versa, starting from the data (D0)-(D6), one can construct a stable
semiabelic variety $V$ together with an ample sheaf $\calO_V(1)$ and
a global section of this sheaf, according to a general construction of
\cite{A2}. Since this construction is crucial to us, we briefly recall
it.

For each cell $\delta$ of the cell decomposition $\Delta$ we consider
$\Cone \delta$, the cone defined by
${1} \times \delta$ in $\RR \times X_{\RR}$.
and to each $\chi=(d,x)\in \XX= \ZZ\oplus X$ in this cone we
associate the sheaf ${\cal M}_{\chi}={\calM}^{\otimes d}\otimes
{c(x)}$, where ${c(x)}$ is the line bundle given by $c(x) \in
\tu{A}$.  The line bundle ${\cal M}_{\chi}$ is rigidified at the
origin of $A$.  Since ${\cal M}_{\chi_1}\otimes{\cal M}_{\chi_2}$ and
${\cal M}_{\chi_1+\chi_2}$ are naturally isomorphic as rigidified
sheaves, we have a semigroup algebra of ${\cal O}_A$-invertible
sheaves labeled by $\Cone\delta$.  The variety ${\tilde V}_{\delta}$
is defined as the $\Proj$ over $A$ of this algebra.  This variety
${\tilde V}_{\delta}$ has a natural ample sheaf ${\cal O}_{{\tilde
    V}_{\delta}}(1)$ and a projection onto $A$.  If two cells
$\delta_1$ and $\delta_2$ intersect along $\delta_{12}$, then the
subvariety ${\tilde V}_{\delta_{12}}$ is contained in both varieties
${\tilde V}_{\delta_1}$ and ${\tilde V}_{\delta_2}$, and we can glue
them along it, together with their sheaves. This will be done in such
a way that the resulting variety is seminormal. In terms of the
$\Proj$ of a ring this can be done as follows. Introduce formal
variables $\zeta_{\chi}$ to label the sheaves ${\cal M}_{\chi}$ that
belong to the cell $\delta$.  Then define a big semigroup algebra $R$,
this time labeled by the union of all the cones
$\Cone \Delta:=\cup_{\delta\in\Delta}\Cone \delta$
(i.e. all $(d,x)$ with $d>0$ and
$(0,0)$), and define the multiplication by setting
$$
\zeta_{\chi_1}\zeta_{\chi_2} = \zeta_{\chi_1+\chi_2}
$$
if $\chi_1$, $\chi_2$ belong to a cone over a common cell of
$\Delta$ (in which case one speaks of cellmates with respect to
$\Delta$), and $0$ otherwise. Let ${\tilde V}= \Proj_A R$. This scheme
carries a natural sheaf ${\cal O}_{\tilde V}(1)$ which restricts to
each ${\tilde V}_{\delta}$ as the sheaf ${\cal O}_{{\tilde
    V}_{\delta}}(1)$.

As we explained above, $\tau$ and $\psi$ define an action of $Y$ on the
algebra $\oplus_{(d,x)\in\ZZ\oplus X} \calM^d\otimes c(x)$.  This
induces a $Y$-action on the algebra $R$, and therefore on $\tilde
V=\Proj_A R$. The variety $(V,{\cal O}_V(1))$ is then defined as the
quotient of $({\tilde V}, {\cal O}_{\tilde V}{(1)})$ by $Y$. This
action is properly discontinuous in the Zariski topology and hence
${\cal O}_V(1)$ is ample due to the Nakai-Moishezon or the Kleiman
criterion.

Finally, the formal power series $\theta=\sum \theta_x\zeta^{(1,x)}$
defines a section of $\calO_V(1)$ as follows. The formal restriction
of $\theta$ to the algebra for the cone $\delta$ defines a section of
$\calO_{V_{\delta}}(1)$. These sections coincide on the intersections
and, by the basic identity in (D5), are compatible with the
$Y$-action.  Therefore, they glue to a global section of $\calO_V(1)$.

For what follows it is important to recall that the variety $V$ has
the following properties (see \cite[1.1]{A2} for more details):
\begin{enumerate}
\item $V$ is seminormal.
\item $V$ is naturally stratified into locally closed strata, and
  there is a 1-to-1 correspondence between the strata of dimension
  $\dim A+i$ and $i$-dimensional cells $\delta$ of the decomposition
  $\Delta$ modulo $Y$.
\item Maximal cells correspond to irreducible components of $V$.
\end{enumerate}

\subsection{Induced degeneration data}
\label{subsec:Induced degeneration data}

Going back to the general setup of the previous subsection, let us now
assume that we have, in addition, a semiabelian variety $G'/S$ with
abelian generic fiber $G'_{\eta}$, and an injective homomorphism
$\varphi:G'\to G$. In this case, the uniformization data for the
variety $(G'_{\eta},\calL'_{\eta}=\varphi^*_{\eta}(\calL_{\eta}))$ can
be read off directly from those of the variety $G_{\eta}$. We would
like to write this down explicitly, since this is exactly the
situation which we encounter with Jacobians and Prym varieties.

By the functoriality of the varieties $\wG,\tu\wG$, we get
homomorphisms $\wG'\to\wG$, $\tu\wG\to\tu\wG'$ which are encoded by
the following commutative diagrams of group schemes over $S$:
    \begin{displaymath}
      \diagram
      X\dto^{c}\rto^{\hat\varphi_T}
      & X'\dto^{c'}
      & Y'\dto^{\tu c'}\rto^{\varphi_Y}
      & Y\dto^{\tu c}
       \\
      \tu A \rto^{\hat\varphi_A}
      & \tu A'
      & A' \rto^{\varphi_A}
      & A.
       \enddiagram
    \end{displaymath}
The homomorphism $T'\to T$ on the torus parts is injective, therefore
$X\to X'$ is surjective. One has $Y'=Y\cap \wG'$, hence $Y'\to Y$ is
injective. The bihomomorphism $\tau:Y\times X\to (\tu c\times
c)^*\calP_{\eta}\inv$ encodes the embedding $Y\to \wG_{\eta}$, and the
same is true for $\tau'$. Therefore, $\tau'$ is merely the ``restriction''
of $\tau$ in the following sense:
$$
\tau'(y',x') = \tau(y',x) \qquad \text{for any } x\mapsto x'. $$
This statement
makes sense for the following reason:
$\tau(y',x)$ is an element of
\linebreak
$c(x)\inv(\tu c(y')) = c(x)\inv\left(\varphi_A(\tu c'(y'))\right)$,
and similarly $\tau'(y',x')$ is an element of
\linebreak
$\varphi_A^*(c(x)\inv) \left(\tu c'(y')\right)$. These are canonically
isomorphic.

The map $\phi':Y'\to X'$ is the composition $Y'\to
Y\overset{\phi}{\to} X\to X'$ and it is injective. Indeed, if
$\phi'(y')=0$ then $\tau'(y',\phi(y'))=1$ by bilinearity. On the
other hand, $\tau'(y',\phi'(y'))=\tau(y',\phi(y'))$, and for
$y'\ne0$ we get a contradiction to the positivity condition on
$\tau$.

The sheaf $\calM'$ is $\varphi_A^*(\calM)$, and we automatically have
$c'\otimes \phi' = \lambda'\otimes\tu c'$, where
$\lambda'=\lambda(\calM'):A'\to\tu A'$. The cubical homomorphism
$\psi'$ is obtained by restriction: $\psi'(y')=\psi(y')$. Again, both
maps take values in canonically isomorphic $\GG_m$-torsors.

On the central fiber, the data (D0)-(D2), and in the case of a
one parameter degeneration also (D3) and (D4), are obtained by
straightforward restriction.

Finally, given a Fourier series $\theta=\sum_{x\in X} \theta_x$ with
$\theta_x\in \Gamma(G_{\eta},\calM_{\eta})$, converging in the
$I$-adic topology, the restriction to $G'_{\eta}$ is given by
$$
\theta'= \sum_{x'\in X'} \theta'_{x'} = \sum_{x'\in X'}
\sum_{x\in X, \, x\mapsto x'} \varphi_A^*(\theta_x).
$$
In this expression, each $\theta'_{x'}$ is itself a power series
converging in the $I$-adic topology, representing an element of
$\Gamma(A'_{\eta},\calM'_{\eta}\otimes c'(x'))$. For as long as
$\Theta_{\eta}$ does not contain $G'_{\eta}$, one has $\theta'\ne0$,
and this can always be achieved by choosing a translated divisor.

\begin{lemma}\label{lem:induced_height}
  For a generic choice of a divisor $\Theta_{\eta}$ in a given
  polarization class the induced height function is given by
  $ H'(x') = \min_{x\mapsto x'} H(x). $
\end{lemma}
\begin{Proof}
  Because of the identity
  connecting $\theta'_{x'}$ and $\theta'_{x'+\phi'(y')}$
  we only have to check this statement for finitely many
  representatives of elements in $X'/Y'$. For each
  $x'$ there are only finitely many $x$ mapping to it with the minimal
  height $H(x)$. We are free to change the cocycle $\psi_0$ by a
  linear function on $Y$, and we are free to multiply the $\theta$'s in
  different
  classes modulo $Y$ by different constants. Hence, it is easy to
  arrange that the sections $\varphi_A^*(\theta_x) $ of the minimal
  height do not cancel each other. \hfill
\end{Proof}
Note that in general we only have an inequality $H'(x')\ge
\min_{x\mapsto x'} H(x)$.

\section{Degenerations of Jacobians and principally polarized Prym varieties}
\newcommand{\oMg}{\overline{M}_g}
\subsection{Degenerations of Jacobians}
\label{subsec:Degenerations of Jacobians}

Let $C_0$ be a stable curve of genus $g$ and let $\hat U$ be the base of the
universal deformation. Then we have a curve $C/{\hat U}$ and a semiabelian
family of Jacobians $JC/{\hat U}$, $JC=\Pic^0_{C/{\hat U}}$, with an abelian
generic fiber. After a finite base change, one can choose an
ample sheaf $\calL$ on $JC$ giving a principal polarization on the
generic fiber. We can then apply
the general construction outlined above to this family.
First we obtain the data (J0)-(J2) which do not depend on the choice of
a one-parameter family.
The data (J3)-(J6)  are \emph{a priori} defined only for a
particular choice of a one-parameter degeneration. However, for the
Jacobians the limit variety
is independent of such a
choice, and so the data (J3)-(J6) depend only on the curve $C_0$
itself, see \cite{A1}.  (This will also follow from Theorem
\ref{thm:indeterminacy_locus} below).  The resulting data are as
follows:

\begin{enumerate}
\item [(J0)] The abelian variety $A_0=JN_0$, where $N_0=\coprod
  N_i$ is the normalization of the nodal curve $C_0$. This is a
  principally polarized abelian variety and the sheaf $\calM_0$ in
  this polarization class can be chosen arbitrarily. Alternatively, on
  the $A_0$-torsor $\prod\Pic^{g_i-1}N_i$ this sheaf is defined
  canonically.
\item [(J1)] The lattice $X=H_1(\Gamma, \ZZ)$, where $\Gamma$ is the
  dual graph of $C_0$, and the same semiabelian variety $G_0=JC_0$
  as in the first section, fitting into an exact sequence
  $$
  1\rightarrow H^1(\Gamma, k^{\ast})\rightarrow G_0=J
  C\rightarrow J N\rightarrow 0,
  $$
  plus the dual semiabelian variety ${^t\!G_0}=G_0$.
\item [(J2)]
  The inclusion $\phi:Y\rightarrow X$ is the identity and
  ${^t\!c}:Y\rightarrow A_0$ is
  defined as $\lambda^{-1}_{A_0}\circ c\circ \phi$.

\item[(J3)] The equivalence class of the bihomomorphism $\tau_0:Y\times
  X \to (\tu c\times c)^*\calP\inv$ is given by the so called Deligne
  symbol, a generalization of the double ratio (see \cite{BM}).
\end{enumerate}

The bilinear function $B:X\times X\rightarrow \ZZ$ is determined by the
monodromy, which is given by the Picard-Lefschetz formula. Explicitly,
we have an embedding
$$
H_1(\Gamma)\subset
C_1(\Gamma)=\oplus_{\operatorname{edges}}\ZZ
e_j.
$$
Let $z_j$, $j=1\ldots q$ be the coordinate functions on
$C_1(\Gamma)$.
On $\hat U$  there
are $q$ transversally intersecting divisors $D_j$, each corresponding
to a node of $C_0$. A one-parameter family $C/S$ gives rise to a map
$S\to\hat U$ and this determines valuations $\alpha_j>0$ along the divisors
$D_j$. The Picard-Lefschetz formula says that the monodromy corresponds
to the form $B=\sum \alpha_j z_j^2$. The positive integers $\alpha_j$ are not
determined by the central fiber $C_0$ alone, but also depend on the
degenerating family.

\begin{enumerate}
\item[(J6)] The cell decomposition $\Delta$ defined
  by $B$ does, however, not depend on the constants $\alpha_j$,
  which are determined
  by the choice of the $1$-parameter family. It is
  obtained by taking the intersections of the standard cubes in $C
  _1(\Gamma, \RR)$ with the subspace $H_1(\Gamma, \RR)$. This follows
  from the simple combinatorial lemma \ref{lem:dicing}
  below and the well-known fact that
  the functions $z_j$ restricted to $H_1(\Gamma,\ZZ)$ form a
  unimodular system of vectors. The corresponding statement in graph theory
  says the following: If $z_{i_1}, \ldots, z_{i_n}$ are linearly independent
  then $\Gamma \setminus \{e_{i_1}, \ldots, e_{i_n}\}$ is a spanning tree.
  This spanning tree defines a basis of $H_1(\Gamma,\ZZ)$ with the
  property $z_{i_s}(e_{i_t})=\delta_{i_si_t}$.
  Therefore $z_{i_1}, \ldots, z_{i_n}$ span $H_1(\Gamma,\ZZ)^*$ over
  the integers $\ZZ$.
\end{enumerate}

Since we are in the principally polarized case the rest of the data
have an auxiliary character, as we have explained in the remarks at the
end of
subsection~\ref{sub:General theory of degenerations of abelian varieties}.

\begin{lemma}\label{lem:dicing}
  Let $X$ be a lattice and $l_j:X\rightarrow\ZZ$ be linear functions
  such that the quadratic form $\sum l^2_j$ is positive definite. Then
  the Delaunay decomposition for the quadratic form $\sum \alpha_j l^2_j$
  is independent of the choice of positive constants $\alpha_j$ if and only
  if the forms $l_j$ define a {\em dicing}: the $0$-skeleton of the
  cell decomposition given by the hyperplanes $\{l_j(x)=n\in\ZZ\}$
  coincides with the lattice $X$.
\end{lemma}
\begin{Proof}
  It is easy to see that the Delaunay decomposition for the quadratic
  form $l^2_{j_0}+\varepsilon\sum_{j\neq j_0} l^2_j,\
  0\le\varepsilon\ll 1$, or, what is the same, for the form
  $Nl^2_{j_0}+\sum_{j\neq j_0} l^2_j, N\gg0$ is a subdecomposition of
  the Delaunay decomposition for the form $l^2_{j_0}$. The latter is,
  of course, formed by the hyperplanes $\{l_{j_0}=n; n\in \ZZ\}$.
  Therefore, if the Delaunay decomposition does not depend on the
  choice of constants $\alpha_j>0$, then it must be a refinement of the
  intersections of dice $\{l_j=n_j; n_j\in\ZZ\}$. Hence this
  intersection cannot have extra vertices, i.e. must be a dicing.

  Vice versa, assume that the linear forms $l_j$ define a dicing. We
  claim that it is a Delaunay decomposition for the quadratic form
  $q=\sum \alpha_j l_j^2$ for any $\alpha_j>0$. Indeed, let $\delta$ be the
  lattice polytope which is the intersection of sets $\{n_j\le l_j\le
  n_j+1\}$ for some $n_j\in \ZZ$. Then for every vertex of $\delta$
  one has $l_j(P)=n_j$ or $n_j+1$ and for any lattice point $Q$, which
  is not a vertex of $\delta$, at least one function $l_j$ has a
  different value. This means that the vertices of $\delta$ all lie on
  an ''empty ellipsoid''
$$
\sum \alpha_j((l_j-n_j-1/2)^2-(1/2)^2)=0
$$
and that all other lattice points lie outside of this ellipsoid.
This is Delaunay's original definition of a Delaunay cell for the form
$q$.  \hfill\end{Proof}

The data above define the unique ``canonical compactified Jacobian''
\linebreak $(\overline{J}_{g-1}C, \Theta_{g-1})$ of \cite{A1}, which is a
stable semiabelic pair. Moreover, this pair defines a point in the
second Voronoi compactification $A_{g}^{\operatorname{Vor}}$ of $A_g$.
If $C_0$ is a stable curve, then this
is the image of $C_0$ under the Torelli map.

\subsection{Degenerations of principally polarized Pryms}
\label{subsec:Degenerations of principally polarized Pryms}

Let $C_0$ be a stable curve with an involution $\iota_0$.
The involution $\iota_0$ gives an involution
on the universal deformation space $\hat U$ of $C_0$, and the fixed
locus $\hat W$ of this involution, which is  regular itself,
represents the deformation of $(C_0,\iota_0)$ together with an
involution. We have a family
$(C,\iota)/\hat W$ with central fiber $(C_0,\iota_0)$ and
a smooth generic fiber $C_{\eta}$. Over $\hat W$ there are two
semiabelian families with abelian generic fibers, namely $JC$ and the subfamily
$P(C,\iota)$. The first one is principally polarized, and that
gives the induced polarization on the second family.

We shall apply subsection~\ref{subsec:Induced degeneration
data} to the subfamily $P(C,\iota)$. This gives us first of all
degeneration data (p0)-(p2), and for the central
fiber the corresponding data (P0)-(P2), which we shall now describe.
We shall see that we have encountered most of these data already in
the first section.

Taking the connected component of the identity of $\ker(1+\iota)$ and
specializing to the central fiber commute, hence the semiabelian
variety $\wG^-_0$ in (P1) is the same as the Prym variety
$P_0$ of the first section (there denoted by $P$), and the abelian variety
$A_0^-$ in (P0) is its abelian part.
The latter was simply denoted by $A$ in the first section.
The ample sheaf $\calM_0^-$ is the
pullback of $\calM_0$ from the Jacobian $JC_0$.
As in section~1 we denote by $[X]^+$, resp.
$[X]^-$ the sublattices of $X=H_1(\Gamma,\ZZ )$ corresponding to the $(\pm
1)$-eigenspaces. Then $\wG^-_0$ is given by the classifying homomorphism
$c^-:X^-\to A_0^-$ where $X^-$ is $X/[X]^+$. This is the classifying
homomorphism which we described in lemma \ref{propclasshom}.

There appears a second lattice in (P1), which we will denote $Y^-$.
It is given by $Y^-=Y\cap \wG^-$ in $\wG$ and we have not yet encountered
this in the first section.
The lattice $Y^-$ contains $2X^-$, and
we have an injective morphism $Y^-\to [X]^-=[H_1(\Gamma,\ZZ)]^-$ which,
however, need not be an isomorphism. In the complex analytic case
$Y^-$ is the image of the natural homomorphism $[H_1(C_t,\ZZ)]^- \to
X=H_1(\Gamma,\ZZ)$, where $(C_t,\iota_t)$ is a topological smoothing of
$(C_0,\iota_0)$. In particular, $Y^-$ can be computed explicitly
in any given concrete example (see also section~5). For curves
$(C_0,\iota_0)$ which correspond to points in $\oRg$ we have $Y^-=2X^-$
for degree reasons.

In (P2), the inclusion $\phi^-:Y^-\to X^-$ is the composition
$Y^-\to [X]^-\to X^-$ and, as in the general case, we have
$c^-\circ \phi^-=\lambda^- \circ \tu c^-$.

The rest of the data (P4)-(P6) require looking at one-parameter
degenerations. To get any kind of uniqueness we must look at the
situation when the induced polarization on the Prym variety is twice
the principal polarization, so that we can analyze the degeneration of
principally polarized Pryms.  We shall, therefore,
for the rest of this section assume that $(C_0,\iota_0)$ is a
nodal curve with an involution which is smoothable to a nonsingular
curve with a fixed-point-free involution. It is easy to show that in
this case $C_0$ has no swapping nodes in the terminology
of the first section, and that the only
points fixed by the involution are the branchwise fixed nodes. If
the genus of the quotient curve $C_0'$ is $g$ then the genus of $C_0$
must be $2g-1$.  We will denote by $\oRg$ the coarse moduli space of
admissible covers, or equivalently of pairs $(C,\iota)$ as above where
$C$ is assumed to be stable. It is not hard to show that
$\oRg$ exists, that it is a compactification of the moduli space of the
moduli space $\calR_{g}$ of smooth genus $g$ curves with an
irreducible \'etale degree-two cover, and that it is proper and
normal.

So, let us fix a one-parameter family $(C,\iota)/S$ so that on the
generic fiber we have a smooth curve with base point free involution.
Our general theory gives us the degeneration data (p0)-(p4) for the
abelian variety $G^-_{\eta}=P(C_{\eta},\iota_{\eta})$ with the
polarization $\lambda^-_{\eta}:G^-_{\eta}\to \tu G^-_{\eta}$
and correspondingly degeneration data (P0)-(P4) for the central fiber.
In our situation these are degeneration data for a twice principally
polarized abelian variety, whereas we are interested in the degeneration
data of the principally polarized abelian variety
whose polarization is given by $\lambda^-_{\eta}=2[\lambda/2]^-_{\eta}$.
We shall denote the latter data by (pp),
resp. (PP), standing for principally polarized Prym.

We can approach this problem backwards. Assuming that we have data
for an ample sheaf $\calL_{\eta}$ on $G_{\eta}$, what are the data for
the sheaf $\calL_{\eta}^2$?  We still have the uniformization of
$G_{\eta}$ as $\wG_{\eta}/Y$. So, the scheme $\wG$ and the embedding
$Y\to \wG$ do not change. Therefore, the schemes $A,\tu A$,
homomorphisms $c$, $\tu c$, and the bihomomorphism $\tau$ do not change.
It is only the cocycle for the sheaf $\calL$ which is doubled. Which
means $\phi=2[\phi/2]$ and $\lambda_A= 2 [\lambda/2]_A$. Therefore, we
get the following data for the central fiber of the principally polarized
Prym:

\begin{enumerate}
\item [(PP0)] The same abelian variety $A_0^-$ as before. The polarization
  $[\lambda/2]^-_0:A_0^-\to \tu A_0^-$ is the unique solution to
  $\lambda^-_0 =2[\lambda/2]^-_0$.
\item [(PP1)] The same lattices $Y^-$ and $X^-$ and homomorphisms
  $c^-:X^-\to \tu A_0^-$ and $\tu c_0^-:Y^-\to A_0$.
\item [(PP2)] $[\phi/2]^-= (1/2)\phi^-:Y^-\to X^-$. This is now an
  isomorphism.
\item [(PP3)] The same $\tau_0^-$ as in the non-principally polarized case.
\item [(PP4)] $[\psi/2]^-_0$ is a solution of the equation ${\psi^-}_0(y)=
  ([\psi/2]^-_0(y))^2$. (This is not unique, and may exist only after a
  finite field extension if the residue field $k$ is not
  algebraically closed.)
\item [(PP5)] Since the polarization is now principal, the section
  is es\-sen\-tial\-ly unique and plays an auxiliary role.
\item [(PP6)] The cell decomposition, which we will denote by $\Delta^-$,
  is the Delaunay decomposition for the bilinear form
  $[B/2]^-:Y^-\times X^-\to \ZZ$, \linebreak $[B/2]^-(y,x)=B^-(y,x)$.
\end{enumerate}

\newcommand{\apgm}{\overline{\operatorname{AP}}_{g-1}}

We have now all the tools to approach the question of the
uniqueness of the limit of the
principally polarized Pryms. The Prym map $\calR_g\to A_{g-1}$,
associating to a smooth curve with a base-point-free involution its
principally polarized Prym variety, extends to a rational map from
$\oRg$ to any compactification of $A_{g-1}$. Of the infinitely many
toroidal compactifications one has a particular functorial importance, namely
the second Voronoi compactification $\Avor_{g-1}$. By \cite{A2} it is
the closure of $A_{g-1}$ in the moduli space $\apgm$ of principally
polarized stable semiabelic pairs, and the latter space is projective.
In \cite{FS} Friedman and Smith have given a series of examples of
pairs $[(C_0,\iota_0)]\in \oRg$ near which the Prym map does not
extend to a regular morphism to $\Avor_{g-1}$ (or any other
reasonable toroidal compactification). (The preprint version
of \cite{FS} contains more details than the published version.) In the
theorem below, we give the complete answer to the extendibility
question. The answer we give is in terms of a
combinatorial condition on the dual graph $(\Gamma,\iota)$ of the
curve $C_0$, together with an involution.

Recall that we have $X^-\otimes\RR\subset X\otimes\RR\subset
C_1(\Gamma, \RR)$. As before, for an edge $e_j$, let $z_j$ be the
corresponding coordinate function on $C_1(\Gamma,\RR)$. We have three
possibilities:
\begin{enumerate}
\item[1.]  $z_j$ is identically zero on $X^-$.
\item[2.]  $z_j: X^-\twoheadrightarrow\ZZ$. This happens if for every
  simple oriented cycle $w \in H_1(\Gamma,\ZZ)$ in the graph $\Gamma$
  $\operatorname{mult}_{e_j}w=1$ implies
  $\operatorname{mult}_{e_{\iota(j)}}w=-1$.
\item[3.]  $z_j:X^-\twoheadrightarrow \frac 1 2\ZZ$. This happens if
  there exists a simple oriented cycle $w \in H_1(\Gamma,\ZZ)$
  with $\operatorname{mult}_{e_j}w=1$ but
                $\operatorname{mult}_{e_{\iota(j)}}w=0$.
\end{enumerate}
The first case is immaterial for us. Define $m_j$ to be 1 or 2
respectively in the second or third case so that we have an
epimorphism from $X^-$ to $\ZZ$.
The combinatorial condition is

\begin{itemize}
\item[(*)] The linear functions $m_j z_j$ define a dicing of the
  lattice $X^-$.
\end{itemize}

Geometrically, the condition $(*)$ can be described as follows.
Consider all translations of the hyperplanes $H_j=\{z_j=0\}$ in
$C_1(\Gamma, \RR)$ through points in $X^-$ and take the
intersection with $X^-_{\RR}$. This defines a cell decomposition and
$(*)$ is fulfilled if the vertices of this cell decomposition are
exactly the points of the lattice $X^-$. In concrete examples this
condition can be checked easily (see section~5 for a number of worked out
examples). Vologodsky \cite{V} has shown that a curve $(C,\iota)$ fails
to fulfill condition $(*)$ if and only if it is the limit of curves of
Friedman-Smith type with at least $4$ nodes.
Here Friedman-Smith type means that $C$ has two
smooth irreducible components $C_1$ and $C_2$ which are fixed under
the involution $\iota$ and that, moreover, the nodes $C_1 \cap C_2$
of $C$ are pairwise interchanged. We are now ready to formulate the main
application of our theory.

\begin{theorem}\label{thm:indeterminacy_locus}
  The following are equivalent:
  \begin{enumerate}
  \item In a neighborhood of the point $0=[(C_0,\iota_0)]$ the
     Prym map $\calR_g\to A_{g-1}$ extends to a morphism from
    $\oRg$ to the second Voronoi toroidal compactification
    $\Avor_{g-1}$.
  \item The limit variety $(\overline{P}_0,\Theta_0)$ depends only on the pair
    $(C_0,\iota_0)$ and not on the choice of a one-parameter
    degeneration.
  \item The cell decomposition $\Delta^-$ depends only on the pair
    $(C_0,\iota_0)$ and not on the choice of a one-parameter
    degeneration.
  \item Condition (*) holds. In this case, the decomposition $\Delta^-$
    coincides with the dicing.
  \end{enumerate}
\end{theorem}
\begin{Proof}
  By the properness of $\apgm$, for every morphism $(S,0)\to (\oRg,0)$
  from a regular one-dimensional $S$ with $(S\setminus 0)\to
  \calR_{g}$, there exists a unique morphism $S\to\Avor_{g-1}$.
  Therefore, (1) implies (2), and clearly (2) implies (3). The
  equivalence of (3) and (4) follows by Lemma~\ref{lem:dicing}. So, we
  only have to show that (3) implies (1).

  Let $Z\subset \Avor_{g-1}$ be the ``image'' of the point $[(C_0,\iota_{0})]$
  under the
  rational map $\oRg\to \Avor_{g-1}$. In other words, if $p:T\to
  \oRg$, $q:T\to\Avor_{g-1}$ is a resolution of this rational map,
  with $p$ proper and birational and $q$ proper, then $Z=q(p\inv(0))$.
  $Z$ is closed, and points of $Z$ correspond to all possible limits
  of Pryms for all one-parameter families $(S,0)\to (\oRg,0)$.  Since
  $\oRg$ is normal, the Zariski Main Theorem implies that $Z$ is
  connected and that the rational map in question is regular in a
  neighborhood of $[(C_0,\iota_0)]$ if and only if $Z$ is a point.

  The structure of any toroidal compactification of $A_{g-1}$ is as
  follows. There is a natural stratification such that each stratum is
  fibered over
  the moduli of the data (D0)-(D2), and each fiber is a quotient of a
  torus by a finite group.  In particular, it is affine. In addition,
  the strata of the second Voronoi compactification correspond to
  Delaunay decompositions modulo $\operatorname{GL}(g-1,\ZZ)$. Since
  the data (D0)-(D2) depend only on $(C_0,\iota_0)$ the image of $Z$
  lies in the union of affine sets as above. If the cell
  decomposition $\Delta^-$ is unique, then there is only one such affine
  set. Since $Z$ is proper and connected it must then be a point and hence
  the Prym map is regular at the point $0=[(C_0,\iota_0)]$.
  \hfill \end{Proof}

\begin{remark}
In order to make the Prym map regular, it is necessary to blow up
$\oRg$. This means that we have infinitely many possible limits of
Prym at the points of indeterminacy. For a description of all possible limits
in the simplest case see the example in \ref{F-S}.
\end{remark}

\section{Degeneration of Pryms with the induced polarization}

Restricting the theta divisor of the Jacobian to the Prym variety defines a
polarization which is usually not a principal polarization.
In this section we shall study possible limits of
Prym varieties with the induced, i.e. in most cases non-principal,
polarization.
Since we now work with the induced polarization we should expect some
relationship with the compactified Jacobian. This we will investigate
in the second part of this section.

\subsection{Construction of the ``middle'' compactified Prym}

Each family of nodal curves with an involution gives us induced data
(p0)-(p4). The problem is that there is no canonical choice of a
section $\theta_{\eta}$ and correspondingly we have no canonical choice
for a cell decomposition for the central fiber. One obtains a
natural candidate for such a cell decomposition by taking the cell
decomposition $\Delta$ for the Jacobian which is given by (J6) and then
intersecting it with $X^-_{\RR}$. As we shall see, this leads indeed to
a degeneration of smooth Prym varieties.
In view of the construction of the cell
decomposition, which is obtained by slicing with $X^-_{\RR}$
(one could also consider slicing with translates of this subspace), we shall
refer to this as the ``middle'' Prym.

Since degenerating families will not appear in the rest of this
section, we will suppress the subscript $_0$ from now on.
Our first aim is to describe data (PM0)-(PM6) which define this
middle Prym. Throughout we will use the letter ``m'' in our
notation in order to indicate ``middle''.

For (PM0)-(PM2) we take the data (P0)-(P2) which, we recall, are universal,
and depend only on the pair $(C,\iota)$. We have already explained our choice
of decomposition, namely we choose

\begin{enumerate}
\item [(PM6)] Define $[\Delta]^-$ to be the cell decomposition
  obtained by intersecting $\Delta$ with the subspace $X^-_{\RR}$.
\end{enumerate}

Clearly, the cells of $[\Delta]^-$ are intersections with $X^-_{\RR}$
of cells $\delta$ of $\Delta$ such that $-\iota(\delta)=\delta$. As a
consequence, the 0-skeleton of $[\Delta]^-$ is contained in the
lattice $X^-$.  In addition, $[\Delta]^-$ is invariant with respect to
the lattice $[X]^-=H_1(\Gamma, \ZZ)\cap H_1(\Gamma, \RR)^-$. Note that
$Y^-\subset [X]^-$ and hence $[\Delta]^-$ is also $Y^-$-invariant.

\begin{enumerate}
\item [(PM3)] For the bihomomorphism $\tau^-_m:Y^-\times X^-\to(\tu
  c^-\times c^-)^*\calP_{A^-\times\tu A^-}\inv$ we take the restriction
  of the bihomomorphism $\tau:Y\times X\to(\tu c\times
  c)^*\calP_{A\times\tu A}\inv$ from the datum (J3) for the Jacobian
  $\overline{J}C$.
\item [(PM4)] We also define the cubical morphism $\psi^-_m$ to be the
  restriction of the morphism $\psi$ in (J4).
\end{enumerate}

It is important to note that, since $[\Delta]^-$ is a subdivision of $\Delta$,
the $\tau$'s and $\psi$'s in the same equivalence class restrict to
equivalent $\tau^-_m$'s and $\psi^-_m$'s.

We will not choose a particular divisor on our
middle Prym, so the
datum (PM5) is not needed.

The above data define a stable semiabelic
variety $[\overline{P}]$ together with an ample
sheaf $[\calL]^-$ according to the
construction we reviewed in subsection~\ref{sub:General theory of
degenerations
of abelian varieties}. It is not a priori clear that this a degeneration of
Prym varieties. This, however, follows from

\begin{theorem}\label{theoremembedding}
  Let $({\cal C},i)/S$ be any one-parameter family of curves with
  involution whose generic fiber is smooth and with special fiber
  isomorphic to $(C,\iota)$. Then, possibly after a finite base change
  $S'\to S$, one can choose a divisor $\Theta_{\eta}$ on the generic
  fiber in the same polarization class as induced by the Jacobian, so
  that the limit stable semiabelic variety of the family $(\cal
  C_{\eta},\Theta_{\eta})$ is the middle Prym $[\overline{P}]$
  defined by $({\cal C},i)/S$.
\end{theorem}
\begin{Proof}
  We first look at the degeneration of Jacobians.  We are free to
  adjust the function $A=\operatorname{val}_t\psi$ by any linear
  function on $Y=X$ without changing the polarization class. After
  a finite base change of degree $2$ (whose effect is multiplying $A$
  and $B$ by 2), we can assume that
  $A^-(y)=B^-(\phi^-(y),y)/2$. Now consider the
  convex hull of the points $(x,H(x))$, $x\in X$.
  It is symmetric with respect to the involution $-\iota$.
If we intersect the
  lower envelope of this hull with the vector space $\RR\oplus
  X^-_{\RR}$, we will get a height function $H''$ on $X^-$ which gives
  the decomposition $[\Delta]^-$. By
  Lemma~\ref{lem:induced_height} it will coincide with an induced
  height $H'$ for a generic choice of the equation $\theta$.
\hfill\end{Proof}

In section \ref{sub:General theory of degenerations of abelian varieties}
(see remark 5 after (D6)) we had
introduced the notion of semi-Delaunay decompositions. These are the
decompositions which play for degenerations with a non-principal polarization
the role which in the case of principal polarizations
are played by Delaunay decompositions. Our above theorem, therefore, has the
following consequence

\begin{corollary}
  The decomposition $[\Delta]^-$ is semi-Delaunay.
\end{corollary}

If $(C,\iota)$ is a nodal curve which can be smoothed to a curve with an
involution
which has at most $2$ fixed points, then the Prym of the
deformation carries twice a principal polarization.
In this case we have two possible degenerations of Prym varieties,
namely the principally polarized Prym variety
$\overline{P}$ and
the middle Prym variety $[\overline{P}]$ with the induced degeneration. In
good cases the limit $\overline{P}$ will be unique and
isomorphic to $[\overline{P}]$ with the only difference that $[\overline{P}]$
carries twice the polarization of $\overline{P}$.
The following result
characterizes when this happens

\begin{proposition}\label{charcterization}
Let $(C,\iota)$ be a nodal curve which can be smoothed to a curve with an
involution with at most $2$ fixed points. Then the following conditions
are equivalent:
\begin{enumerate}
  \item
The middle Prym variety $[\overline{P}]$ is as a
stable semiabelic variety isomorphic
to $\overline{P}$ and carries twice its polarization
 \item
The following condition is fulfilled:
\begin{itemize}
\item [(**)] $[\Delta]^-$ is a dicing with respect to
  the lattice $2X^-$.
\end{itemize}
\end{enumerate}
\end{proposition}
\begin{Proof}
  By the assumption that $(C,\iota)$ can be smoothed to a curve $C$
  with an involution with at most $2$ fixed points we have the
  equality of lattices $Y^-=[X]^-=2X^-$. Clearly the second condition
  is equivalent to $[\Delta]^-=2\Delta^-$.
  On the other hand,
  following the discussion of conditions (PP0-PP6) in section
  \ref{subsec:Degenerations of principally polarized Pryms}, the
  principally polarized Prym $\overline{P}$ with twice the principal
  polarization is given by the cell decomposition $2{\Delta}^-$ and
  $Y^-=2X^-$. Moreover, by construction the bilinear forms $\tau^-_m$
  and $\tau|_{2X^-\times X^-}$ coincide. Hence, $[\overline{P}]\simeq
  \overline{P}$ if and only if $[\Delta]^-=2\Delta^-$.
% it follows from the
%construction of the central fiber according to the procedure
%outlined at the end of
%section \ref{sub:General theory of degenerations of abelian varieties}
%that this is equivalent to the first condition.
\hfill \end{Proof}

\begin{remarks}
\begin{enumerate}
\item
Obviously condition (**) is stronger than (*). For a discussion
of an example where (*) holds, but condition (**) fails,
see example \ref{3-comp}.
\item For any given nodal curve $(C,\iota)$ it is easy to verify
whether (**) is fulfilled.
\item In the appendix Vologodsky \cite{V} shows that a
stable curve $(C,\iota)$
fails condition (**) if and only if it is a degeneration of a Friedman-Smith
type curves with at least 2 nodes.
\end{enumerate}
\end{remarks}

\subsection{Relation with the compactified Jacobian}

The last question which we would like to address is the relationship
between the
variety $[\overline{P}]$ and $\overline{J}C$, as
motivated by the smooth case.

We start with some remarks concerning maps of toric varieties.
Let $\delta$ be a polytope in $\Delta$ and consider the lattice
polytope $\delta^-=\pi^-(\delta)$ in the lattice $X^-=\pi^-(X)$.
Then we have two associated algebras:
\begin{enumerate}
\item The $\calO_A$-algebra $R=R_{\delta}$ which is the subalgebra in
  $\oplus_{(d,x)\in\ZZ \oplus X} \calM^d\otimes c(x)$ corresponding to
  the lattice points in the cone over the polytope $\delta\subset
  (1,X)$
(here $A=JN$, the Jacobian of the normalization of $C$).
This algebra is graded by the semigroup $\Cone\delta\subset
  \ZZ\oplus X$.
\item The analogous $\calO_{A^-}$-algebra $R^-=R_{\delta^-}$.
\end{enumerate}
Since the morphism $A^-\to A$ is finite (recall from
section \ref{Construction of $P$}
that its kernel is $\im(\ZZ/2\ZZ)^k$, where $k$ is the
number of $(+1)$-blocks in the decomposition of the involution $\iota$
on $X$), $A^-$ is affine over $A$. As an $\calO_A$-algebra, $R^-$ is
graded by the semigroup $\Cone\delta^- \oplus \im(\ZZ/2\ZZ)^k\subset
\ZZ\oplus X^-\oplus\im(\ZZ/2\ZZ)^k$.  The homomorphism $R\to R^-$ is a
homomorphism of graded algebras, and the homomorphism $X\to X^-\oplus
\im(\ZZ/2\ZZ)^k$ comes from equation~(\ref{eq:lattice_for_Prym}) after
Lemma \ref{lemsemiabelian}.  For both gradings
we have an invertible $\calO_A$, resp. $\calO_{\im A^-}$-module
($\im A^- = K_N$ as in the first section) in every homogeneous degree.

\begin{lemma}\label{lem:graded_algs}
\begin{enumerate}
\item $\Spec R^- \to \Spec R$ is a closed embedding if and only if the
  semigroup homomorphism $\Cone\delta \to \Cone\delta^- \oplus
  \im(\ZZ/2\ZZ)^k$ is surjective.
\item $V_{\delta^-}=\Proj R^- \to V_{\delta}=\Proj R$
 is a closed  embedding if
  and only if the image of $\Cone\delta \to
  \Cone\delta^- \oplus \im(\ZZ/2\ZZ)^k$ gives everything in high enough
  degrees $d\ge d_0$.
\item The morphism from the main, i.e. the dense, stratum of
  $V_{\delta^-}$ to the main stratum of $V_{\delta}$ is injective if
  and only if the homomorphism $\RR\delta\cap X \to (\RR\delta^-\cap
  X^-) \oplus \im(\ZZ/2\ZZ)^k$ of abelian groups is surjective.
\item the morphism
  $V_{\delta^-}\to V_{\delta}$ is finite.
\end{enumerate}
\end{lemma}
\begin{Proof}
The first three statements are applications of
basic facts about $\Proj$'s of graded algebras and
toric varieties. The fourth one follows since, in any case, $\delta\to\delta^-$
 is a surjective map
  of polytopes. Hence the saturation of the image of the
corresponding homomorphism of semigroups is everything. Thus
 the algebra $R^-$ is finite over $R$, and the morphism of $\Proj$'s is
well defined and finite. \hfill
\end{Proof}

For every polytope $\delta^-$ in $[\Delta]^-$ there exists a unique minimal
polytope $\delta$ in $\Delta$ such that $\delta^-=\delta\cap X^-_{\RR}$. Under
the involution $-\iota$ the polytope $\delta^-$ and the whole decomposition
$\Delta$ map to themselves. Therefore, the polytope $\delta$ has to be
invariant
under the involution as well. This implies that under $\pi^-:X_{\RR}\to
X^-_{\RR}$ the cell $\delta$ maps surjectively to $\delta^-$, and so we are in
the situation of the above lemma. Putting all polytopes together and applying
the compatible actions of the lattices $Y$ and $Y^-$ now leads to a finite
morphism from $[\overline{P}]$ to $\overline{J}C=\overline{J}_{g-1}C$.
Already in
simple examples, however, it  turns out that this morphism in not optimal.
The basic example is the situation considered by Beauville \cite{B}
(see also Example \ref{Beauville type examples} in the next section)
when the Prym
variety is in fact abelian. Without a shift we get the finite
morphism $P\to K_N$ from section 1, whose degree is some power of 2,
and which need not be an embedding. In
order to discuss the optimal choice for the shift $v$ we start with the
following

\begin{lemma}
  Let $v$ be any vector in $(1/2) (X\cap [C_1(\Gamma,\ZZ)]^+ )$.
  Then the cell decompositions $[\Delta]^- +v$ and $\Delta \cap (X_{\RR}^- +v)$
  coincide.
\end{lemma}
\begin{Proof}
Let $u\in H_1(\Gamma,\RR) \subset C_1(\Gamma,\RR)$ be an arbitrary vector.
Write it as $\sum u_je_j$ using the standard basis in $C_1(\Gamma,\ZZ)$,
labeled by the edges
$e_j$. Let us first determine the unique cell $\delta(u)$ in $\Delta$ which
contains $u$ in its relative interior. After shifting $u$ by an element of
$H_1(\Gamma,\ZZ)$ we can assume that $-1< u_j <1$ for all $j$. Then $u$ lies
in the relative interior of a standard Euclidean cube in $C_1(\Gamma,\RR)$
determined by a system of $|J|$ equalities or inequalities which are: for each
$j$, $z_j=0$ if $u_j=0$, $-1\le z_j\le 0$ if $u_j$ is negative, and $0\le z_j
\le 1$ if $u_j$ is positive.
Obviously, the cell $\delta(u)$ is the intersection of this
cube with $H_1(\Gamma,\RR)$. Denote by $\supp (u)$ the collection of edges
$e_j$
with $u_j\ne0$, and by $\Gamma(u)$ the spanning subgraph of $\Gamma$ which has
these edges.
Then it is easy to see that the cell $\delta(u)$ spans the sublattice
$H_1(\Gamma(u),\ZZ) \subset H_1(\Gamma,\ZZ)$.

Now fix the vector $v$. The two decompositions of our lemma coincide if
and only if for
every $u\in X^-_{\RR}$ the negative parts of $H_1( \Gamma(u))$ and
$H_1\left(\Gamma(u+v)\right)$ coincide.
But this is clear: for every $e_s$ appearing in
$v=\sum v_je_j$ with a non-zero coefficient, one has $\iota(e_s)=e_s$ and
such an edge does not appear in the support of any $u$  with $\iota(u)= -u$.
\hfill
\end{Proof}

The cubical morphism $\psi_m^-$ is well suited for the inclusion
$Y^-\to Y$ but not for the inclusion shifted by $v$. Since $v$ is
half-integral, we either have to make a choice out of $2^{\,\dim X^-}$
possibilities or work with the sheaf $\calL^2$.
{}For the latter purpose, we define

\begin{enumerate}
\item [(PM4)$_2$] the cubical morphism
$\psi^-_{m,2}:Y^-\to (\tu c^-)^*\calM^{*-2}$ is given by the formula
$ \psi^-_{m,2}(y^-) = \psi_m^2(y^-) \tau(y^-,2v).$
\end{enumerate}

Note that this choice only affects the ample sheaf
on $[\overline{P}]$,
but not
the variety itself.

\begin{theorem}\label{theoremmiddlePrym}
  For any vector in $(1/2)( X \cap
  [C_1(\Gamma,\ZZ)]^+ )$ there exists a finite morphism $f_v$ from the
  ``middle'' Prym variety $[\overline{P}]$ to the canonically
  compactified Jacobian $\overline{J}C$.
%  Assume $Y^-=[X]^-$. Then for any vector in $(1/2)( X \cap
%  [C_1(\Gamma,\ZZ)]^+ )$ there exists a finite morphism $f_v$ from the
%  ``middle'' Prym variety $[\overline{P}]$ to the canonically
%  compactified Jacobian $\overline{J}C$. Without the assumption
%  $Y^-=[X]^-$ there is a finite morphism to a $[X]^-/Y^-$-Galois cover
%  of $\overline{J}C$ for each $v$.
\end{theorem}

\begin{Proof}
  Take a cell $\delta^-$ of $[\Delta]^-$ and let $\delta$ be the
  minimal cell in $\Delta$ containing $\delta^-+v$. As before, we see
  that under $\pi^-$ the polytope $\delta$ surjects to $\delta^-$.
  Therefore, we have a finite morphism $V_{\delta^-}\to V_{\delta}$.
  By construction, these morphisms coincide on the intersections of
  cells, and so glue to a finite morphism $\tilde V^-\to \tilde V$.
  The functions $\tau,\psi^2$ define the $Y$-action on $\left(\tilde
    V,\calO(2)\right)$ and, in particular, the action by any subgroup
  $Y'\subset Y$ such that $Y'\cap [X]^-= Y^-$. Similarly, the
  functions $\tau^-_{m},\psi^-_{m,2}$ define the $Y^-$-action on
  $\left(\tilde V, \calO(2)\right)$.  Since these functions were
  chosen to be compatible, we have a finite morphism on the quotients
  $V^-=\tilde V^-/Y^- \to \tilde V/Y'$.  This gives a finite map from
  $ [\overline{P}]$ to a $[X]^-/Y^-$-Galois cover of $\overline{J}C$.
  Dividing by $[X]^-/Y^-$, we get the required finite morphism.
  \hfill
\end{Proof}

Of course, two half-integral
shifts that differ by an integral vector give the same
morphism $f_v$. As we said before,
the basic motivation for making a shift $v$ is the situation
considered in \cite{B} where the Prym
variety is in fact abelian
and where, without a shift, we will in general only obtain a finite map,
but not an embedding.
In this case we obtain an embedding if we shift by a vector
$v=\sum_{\text{all } j}e_j/2$. With this in mind, we give the following

\begin{definition}
   We call a vector in $(1/2) (X \cap [C_1(\Gamma,\ZZ)]^+ )$
a \emph{maximal half-shift} if the support of $v$ consists of all edges
$e_j$ with $\iota(e_j)=e_j$, i.e. if it spans $[C_1(\Gamma,\ZZ)]^+$.
\end{definition}

\begin{lemma}
  Assume that the only fixed points of the involution on $C$ are
branchwise fixed nodes.
Then a maximal half-shift
exists, and it is unique up to an element of $H_1(\Gamma,\ZZ)$.
\end{lemma}
\begin{Proof}
  Under our assumption the fixed points of
the involution on the normalization $N$
are precisely the preimages of the branchwise fixed nodes. On the graph
$\Gamma$ these nodes correspond to edges $e_j$ which are fixed under $\iota$.
Let $\Gamma^+$ be the spanning subgraph of $\Gamma$ which has these nodes.
Since the number of fixed points of an involution on a smooth curve is even,
the degrees of the vertices in $\Gamma^+$ are all even. By a basic fact in
graph theory there exists an oriented eulerian cycle $\sum e_j$, and half of
it gives us a maximal half-shift. Moreover, two such cycles obviously
differ by an even integral vector in $H_1(\Gamma,\ZZ)$, so half of it
will be integral. \hfill
\end{Proof}

As a consequence of condition (3) of Lemma~\ref{lem:graded_algs} we obtain
\begin{proposition}
  Assume the map $X\cap [C_1(\Gamma,\ZZ)]^+ \to
\im(\ZZ/2\ZZ)^k$ to be surjective. Then the
morphism $f_v$ for a maximal half-shift $v$  is generically injective.
\end{proposition}

In general, the question of when the morphism $f_v$, in particular for a
maximal
half-shift, is an embedding, appears to be combinatorially quite involved,
and is deserving of a separate study. We leave this question for  another
place and time.

%%%%%%%%%%%%%%%%%%%%%%%%%
%%%%%%%%%%%%%%%%%%%%%%%%%

\section{Examples}
In this section we illustrate our theory by discussing degenerations of Prym
varieties in several examples. We shall mostly
concentrate on a description of the
combinatorial data, i.e. the lattices $Y^-,X^-$ and $[X]^-$ and the cell
decompositions $\Delta^-$ and $[\Delta]^-$. This gives the combinatorial
structure of the toric part of the degenerate Prym variety.

The examples \ref{Beauville type examples} to \ref{3-comp}
all correspoind to points in
$\oRg$ (at least if there are no smooth fixed ponts).
The examples \ref{Case 3.3} and \ref{exchanged-components} are optimally
behaved in the sense that condition (**), and hence also (*), is fulfilled.
These examples demonstrate different behaviour of the toric parts. Example
\ref{F-S} explains the Friedman-Smith examples and example \ref{3-comp} is the
easiest example in $\oRg$ where (*) holds, but (**) fails. Here we explain
the relationship between the principally polarized and the middle Prym. The
remaining examples treat mixed cases.

\subsection{Beauville type examples}\label{Beauville type examples}
Let $C$ be a nodal curve with an involution $\iota$ which has only
branchwise fixed nodes.  In particular $\iota$ maps every component $C_i$
of $C$ to
itself and the involution $\iota$ acts as the identity on $C_1(\Gamma,\ZZ)$
and hence also
on $H_1(\Gamma,\ZZ)$. In this case $X^-=Y^-=\{0\}$ and the Prym variety is
an abelian variety. If there are no smooth fixed points,
then these are the principal examples treated by Beauville in \cite{B}. The
presence of smooth fixed points does not change $X^-=Y^-=\{0\}$, but if we have
more than 2 smooth
fixed points the polarization on the Prym variety will no longer be twice a
principal polarization.

\subsection{Non-fixed nodes only}
\subsubsection{ Case of irreducible $C$}\label{Case 3.3}
Let $C$ be an irreducible nodal curve with an even number of double points
$Q_1,\ldots,
Q_{2a}$ which are pairwise exchanged by $\iota$, say
$\iota(Q_j)=Q_{2a-j+1}$. The graph
$\Gamma(C)$ has one vertex $v$ and $2 a$ loops $e_1,\ldots,e_{2a}$ around
$v$.\hfill

\vspace*{0.5cm}
\setlength{\unitlength}{1cm}
\begin{figure}[h]
\begin{center}\unitlength 1cm
\epsfig{file=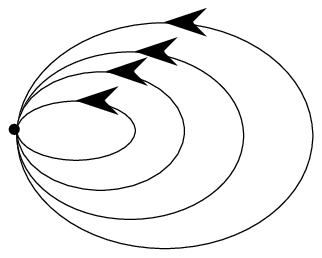}
\put(-3.5,1.0){$v$}
\put(0.1,1.){$e_4$}
\put(-0.7,1.){$e_3$}
\put(-1.3,1.0){$e_2$}
\put(-1.85,1.0){$e_1$}
\end{center}
\end{figure}
\vspace*{0.5cm}
\noindent In this case the involution $\iota$ acts on
$H_1(\Gamma,\ZZ)=C_1(\Gamma,\ZZ)$ by
$\iota(e_j)=e_{2a-j+1}$ and is,  therefore, with respect to the basis $e_1,
e_{2a}, e_2,
e_{2a-1},\ldots,e_{a}, e_{a+1}$ given by the matrix
$$
\begin{pmatrix}
\begin{pmatrix}
0 & 1\\ 1 & 0
\end{pmatrix}&&&\\
&\begin{pmatrix}
0 & 1\\ 1 & 0
\end{pmatrix}&&\\
&&\ddots&\\
&&&\begin{pmatrix}
0 & 1\\ 1 & 0
\end{pmatrix}
\end{pmatrix}.
$$
The decomposition $\Delta$ of $H_1(\Gamma,\ZZ)=C_1(\Gamma,\ZZ)$ is given by the
standard cubes. In this
case we have
$$
[H_1(\Gamma,\ZZ)]^-=[X]^-=\langle e_1-e_{2a},\ldots, e_{a}-e_{a+1}\rangle=2X^-.
$$
The induced decomposition $[\Delta]^-$ consists also of standard cubes, i.e.
the cube
spanned by $e_1-e_{2a},\ldots, e_{a}-e_{a+1}$  its faces
and the
$2X^-$-translates of these cells. It is Delaunay with respect to $2X^-$. This
implies that condition $(**)$, and hence also
condition $(*)$, is fulfilled. Since the smoothing of non-fixed nodes
does not
introduce fixed points on $C_t$, we have already for degree reasons that
 $Y^-=2X^-=[X]^-$. The decomposition $\Delta^-$ consists of cubes with vertices
in $X^-$ and we have $[\Delta]^-=2\Delta^-$.
As semiabelic varieties (without the polarization)
$\overline{P}=[\overline{P}]$ and if there
is no abelian part this variety is isomorphic to
$(\PP^1)^a$ where \lq\lq opposite\rq\rq sides are glued and the gluing is
defined by the form $\tau^-$. If there is an abelian part, then
the normalization of $\overline{P}=[\overline{P}]$ is a $(\PP^1)^a$-bundle
over the abelian part.

\subsubsection{Two exchanged components}\label{exchanged-components}
Assume that $C=C_1\cup C_2$ with $C_1$ and $C_2$ intersecting in
$Q_1,\ldots, Q_{2a}$.
The involution $\iota$ interchanges $C_1$ and $C_2$ and acts on the nodes by
$\iota(Q_j)=Q_{2a+1-j}$. Then the graph $\Gamma(C)$ has 2 vertices and $2
a$ edges and
looks as follows
\setlength{\unitlength}{1cm}
\begin{figure}[h]
\begin{center}
\epsfig{file=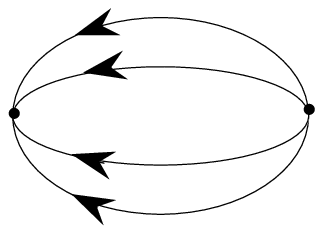}
%E2a.eps
\put(-2,2.3){$e_1$}
\put(-2,1.8){$e_2$}
\put(-1.9,1.0){$\vdots$}
\put(-4,1.0){$v_1$}
\put(-0.1,1.0){$v_2$}
\put(-2,0.4){$e_{2a-1}$}
\put(-2,-0.2){$e_{2a}$}
\end{center}
\end{figure}

\noindent The involution $\iota$ acts by $\iota(e_j)=-e_{2a+1-j}$ and the
following elements
define a basis of $H_1(\Gamma,\ZZ)$:
$$
\begin{array}{l}
h_1 = e_1-e_{2a},\\
h_{2m} = e_m-e_{2a-m};\ m=1,\ldots, a-1\\
h_{2m+1} = e_{m+1}-e_{2a-m+1};\ m=1,\ldots, a-1.
\end{array}
$$
With respect to this basis the involution $\iota$ is given by the matrix
$$
\iota=
\begin{pmatrix}
1&&&\\
&\begin{pmatrix}
0 & 1\\ 1 & 0
\end{pmatrix}&&\\
&&\ddots&\\
&&&\begin{pmatrix}
0 & 1\\ 1 & 0
\end{pmatrix}
\end{pmatrix}.
$$
We find that
$$
[H_1(\Gamma,\ZZ)]^-=[X]^-=\langle h_2-h_3,\ldots, h_{2a-2}-h_{2a-1}\rangle =
2X^-.
$$
If we set
$$
l_1=\frac 1 2 (h_2-h_3), l_2=\frac 1 2(h_4-h_5),\ldots,
l_{a-1}=\frac 1 2(h_{2a-2}-h_{2a-1}),
$$
then
$$
X^-=\langle l_1,\ldots, l_{a-1}\rangle.
$$
Since non-fixed nodes do not contribute to fixed points on $C_t$ we find that
$Y^-=2X^-$ in all of these cases.
%\newpage
In  case $a=1$ we have $X^-=\{0\}$. If $a=2$ we find the following
picture
\newpage
\setlength{\unitlength}{1cm}
\begin{figure}[h]
\begin{center}
\epsfig{file=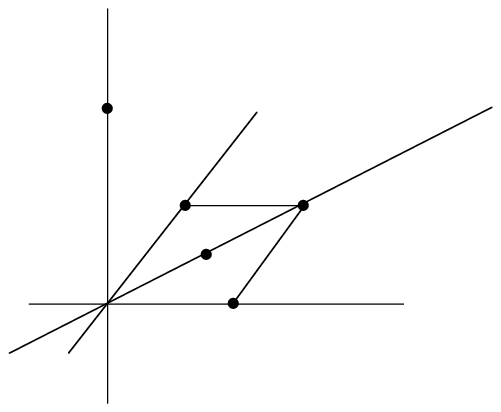}
%E3.1_.eps}
\put(-3.7,3.1){$h_1$}
\put(-3.9,2.2){$-h_3$}
\put(-1.7,1.8){$h_2-h_3$}
\put(-2.7,0.5){$h_2$}
\put(0.5,3){$X^-_{\RR}$}
\end{center}
\end{figure}

\noindent In this case $[\Delta]^-$ is the unique Delaunay decomposition
with respect to
the rank 1 lattice $2X^-$ and similarly $\Delta^-$ is the unique Delaunay
decomposition for $X^-$. As semiabelic varieties $\overline{P}=[\overline{P}]$
and this is a nodal curve (with a polarization of degree $1$, resp. $2$).

Let $a=3$. Then
 $X^-=\langle l_1, l_2\rangle$ and the decomposition $[\Delta]^-$ of
$X^-_{\RR}$
looks as follows (the broken lines indicate $\Delta^-$)

\setlength{\unitlength}{1cm}
\begin{figure}[h]
\begin{center}
\epsfig{file=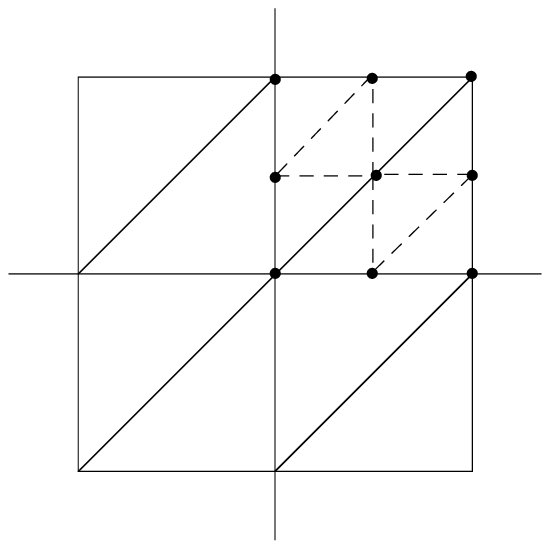}
\put(-2.5,4.8){$2l_2$}
\put(-2.5,3.5){$l_2$}
\put(-1.9,2.2){$l_1$}
\put(-0.5,2.2){$2l_1$}
\put(0.1,4.3){$X^-_{\RR}$}
\end{center}
\end{figure}

\noindent Again $[\Delta]^-=2\Delta^-$ and both decompositions
are Delaunay with respect to
$2X^-$, resp. $X^-$. If the curves $C_1$ and
$C_2$ are rational we obtain $2$ copies of $\PP^2$ each with a polarization of
degree $2$, resp. $1$.
Note that the $2$ projective planes are glued in exactly the same way as in the
case of the
compactified Jacobian of the dollar sign curve of genus $2$.
Here condition (**), and hence also (*), is fulfilled.

For general $a$ let
$y_1,\ldots, y_{a-1}$ be the dual coordinates with respect to $l_1,\ldots,
l_{a-1}$. Then
$[\Delta]^-$ is given by
$$
\begin{array}{ll}
y_i=c\in 2\ZZ; & i=1,\ldots, a-1,\\
y_i-y_{i+1}=d\in 2\ZZ; & i=1,\ldots, a-2.
\end{array}
$$
This gives a Delaunay decomposition for $2X^-$ and, in particular,
condition (**), and thus also (*), is fulfilled.
The building blocks of the toric part are
projective spaces $\PP^{a-1}$.

\subsubsection{The Friedman-Smith examples}\label{F-S}
Let $C=C_1\cup C_2$ consist of two irreducible components intersecting in
an even number
of nodes $Q_1,\ldots, Q_{2a}$. Assume that $C_1$ and $C_2$ are fixed by the
involution
$\iota$ which, however, interchanges the nodes pairwise, i.e.
$\iota(Q_j)=Q_{2a-j+1}$. This
situation can, for example, be realized as follows. Let $C_1$ and $C_2$ be
two elliptic
curves and choose non-zero 2-torsion points $\tau_1,\tau_2$. Choose general
points
$R_1,\ldots, R_a$ on $C_1$ and $S_1,\ldots, S_a$ an $C_2$. Next identify
the points $R_k$
and $S_k$ for $k=1,\ldots, a$ as well as $R_k+\tau_1$ and $S_k+\tau_2$ for
$k=1,\ldots,a$.
(We choose the points $R_k$ such that $R_{k_1}+\tau_1\neq R_{k_2}$ for all
$k_1, k_2$ and
similarly for the points $S_k$). Then the involutions $x\mapsto x+\tau_i;
i=1,2$ on $C_i$
define an involution $\iota$ on $C$ which has the required properties.

As before the graph $\Gamma(C)$ has 2 vertices and $2a$ edges
\setlength{\unitlength}{1cm}
\begin{figure}[h]
\begin{center}
\epsfig{file=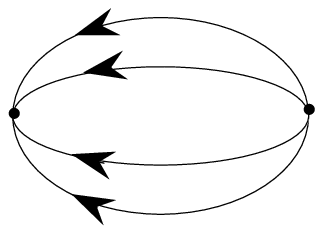}
\put(-2,2.3){$e_1$}
\put(-2,1.8){$e_2$}
\put(-1.9,1.0){$\vdots$}
\put(-4,1.0){$v_1$}
\put(-0.1,1.0){$v_2$}
\put(-2,0.4){$e_{2a-1}$}
\put(-2,-0.2){$e_{2a}$}
\end{center}
\end{figure}

\noindent The difference is that this time $\iota(e_j)=e_{2a-j+1}$. We can
define a basis
of
$H_1(\Gamma, \ZZ)$ by setting
$$
\begin{array}{l}
h_1= e_1-e_{2a},\\
h_{2m}= e_m-e_{2a-m};  m= 1,\ldots, a-1\\
h_{2m+1}= e_{2a-m+1}-e_{m+1},  m=1,\ldots, a-1.
\end{array}
$$
With respect to this basis the involution $\iota$ is given by the matrix
$$
\iota=
\begin{pmatrix}
-1&&&\\
&\begin{pmatrix}
0 & 1\\ 1 & 0
\end{pmatrix}&&\\
&&\ddots&\\
&&&\begin{pmatrix}
0 & 1\\ 1 & 0
\end{pmatrix}
\end{pmatrix}.
$$
To describe the lattices $[X]^-$ and $X^-$ we define
$$
l_0=h_1,\ l_j=\frac 1 2 (h_{2j}-h_{2j+1}); \ j=1,\ldots, {a-1}.
$$
Then
$$
[X]^-=\langle l_0, 2l_1,\ldots,2l_{a-1}\rangle,\quad X^-=\langle l_0,
l_1,\ldots,
l_{a-1}\rangle.
$$
We first consider the case $a=1$. This case is special in the sense that
$X^-=[X]^-$. The
decomposition $[\Delta]^-$ is the unique Delaunay decomposition with
respect to $X^-$. In
particular it is not Delaunay with respect to $2X^-$. On the other hand
condition (*) is
fulfilled since $X^-$ has rank 1.

Let us now consider the case $a=2$. The decomposition $[\Delta]^-$ then
looks as follows

\setlength{\unitlength}{1cm}
\begin{figure}[h]
\begin{center}
\epsfig{file=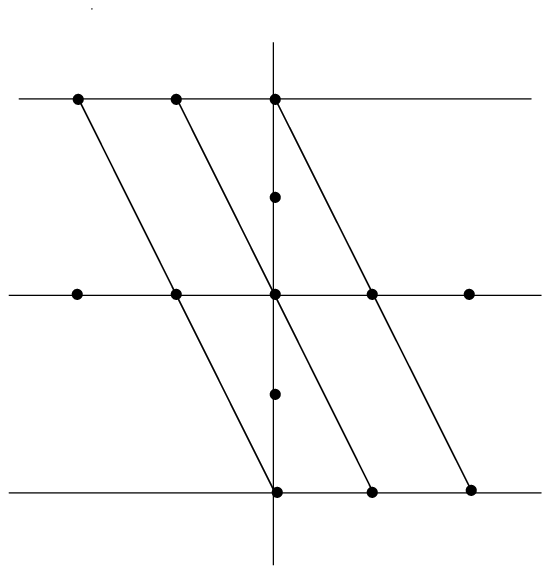}
\put(-1.7,2.85){$l_0$}
\put(-0.8,2.85){$2l_0$}
\put(-2.55,3.5){$l_1$}
\put(-2.55,4.35){$2l_1$}
\end{center}
\end{figure}

\noindent We can see immediately from this picture that condition (*) is not
fulfilled. Similarly for arbitrary integers $a\ge 2$ the decomposition
$[\Delta]^-$ is not
Delaunay with respect to $2X^-$, nor is (*) fulfilled. If we have no smooth
fixed points on $C_1$ and $C_2$, then $Y^-=2X^-$. This is not true in general.
If e.g. $C_1$ and $C_2$ are genus $2$ curves and $\iota$ has $2$ smooth
fixed points on each of the curves $C_i$, then the abelian part of the Prym is
$(2,2)$-polarized and $X^-/Y^-=\ZZ/2$.

In this example one can see easily that there are different possible
limits for the principally polarized Prym variety depending on the
degenerating family. Recall that the form $B$ on the lattice $X=\ZZ e_1 +
\ZZ e_2 + \ZZ e_3 + \ZZ e_4$ is of the form
$B=\sum_{j=1}^4 \alpha_jz^2_j$, where
the $\alpha_j$ depend on the degenerating family.
The induced form on the lattice
$X^-=\ZZ l_1 + \ZZ l_2$ is given by the matrix
$$
B^-=\begin{pmatrix} \alpha_1 + \alpha_4 & (\alpha_1 + \alpha_4)/2
\\ (\alpha_1 + \alpha_4)/2 &
(\alpha_1 + \alpha_2 + \alpha_3 + \alpha_4)/4  \end{pmatrix}.
$$
If $\alpha_j=1$ for $j=1, \ldots, 4$ this leads to $B^- =
\begin{pmatrix} 2 & 1\\ 1 & 1\end{pmatrix}$ which is the sum of two
semi-positive definite matrices and hence the corresponding Delaunay
decomposition of the plane consists of squares.
Geometrically this means the following. Let us, for simplicity, assume
that we have no abelian part (i.e., the curves $C_i$ are
elliptic). Then the normalization of
$\overline{P}$ is a quadric $\PP^1 \times \PP^1$
and $\overline{P}$ is obtained from this quadric by identifying opposite sides
of the square given by the coordinate points via the relation $x \sim bx$,
where $b$ is a non-zero constant. (If there is an abelian part we have a
fibration of quadrics over it.) Incidentally, the constant $b$ corresponds to
the form $\tau$.
On the other hand if
$\alpha_1=\alpha_4=1$ and $\alpha_2=\alpha_3=3$ then $B^- =
\begin{pmatrix} 2 & 1\\ 1 & 2\end{pmatrix}$ which is the sum of three
semi-positive definite forms and the corresponding Delaunay decomposition
of the plane consists of triangles. Geometrically this correspond to two
copies of $\PP^2$ glued as in the Jacobian of a dollar sign curve of genus $2$.
The parameters $b$ and $b^{-1}$ give isomorphic varieties.
Altogether we obtain, as $b$ varies, a $\PP^1$ worth of degenerations,
where $b= 0, \infty$ belong to the case of the union of two planes.

\subsubsection{An example with $3$ components}\label{3-comp}
Consider a nodal curve $C$ which looks as follows
\setlength{\unitlength}{1cm}
\begin{figure}[h]
\begin{center}
\epsfig{file=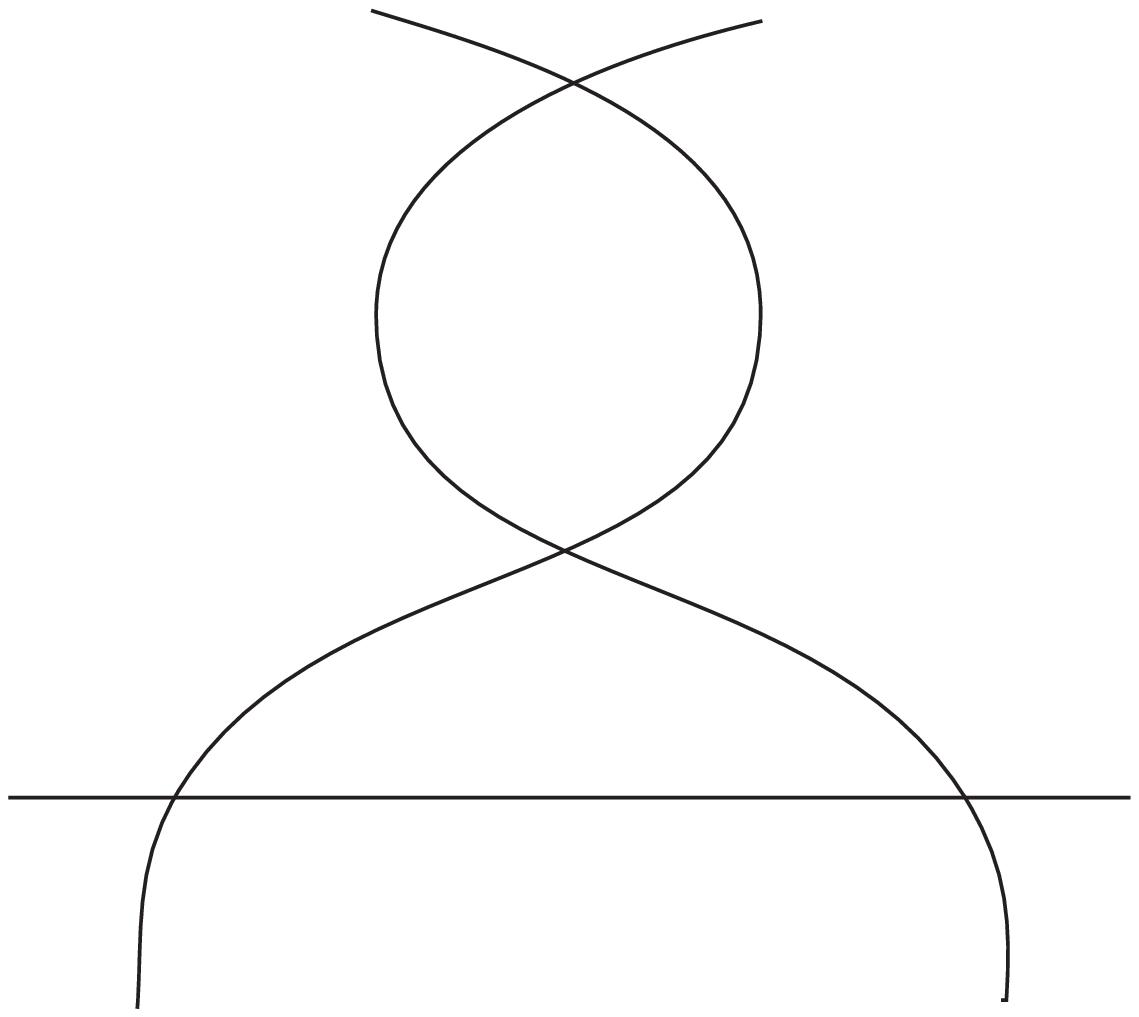,scale=0.5}
\put(-1.9,4.5){$Q_1$}
\put(-1.9,2.4){$Q_2$}
\put(-1.6,5.3){$C_2$}
\put(-4.5,5.3){$C_1$}
\put(-5.6,1.4){$Q_3$}
\put(-0.6,1.4){$Q_4$}
\put(0.5,1.0){$E$}
\end{center}
\end{figure}

\noindent For the involution $\iota$ we want to assume that $\iota(C_1)=C_2,
\iota(Q_1)=Q_2$ and
$\iota(Q_3)=Q_4$. Then we must have $\iota(E)=E$.
This is easy to obtain if we take e.g. $E$ to be an
elliptic curve and
$\iota$ a fixed point free involution on $E$ which interchanges the points
$Q_3$ and
$Q_4$. In this case $\iota$ has no smooth fixed points. We orient the graph
$\Gamma(C)$ of $C$ as follows
%\newpage

\setcounter{totalnumber}{5}
\setlength{\unitlength}{1cm}
\begin{figure}[h]
\begin{center}
\includegraphics{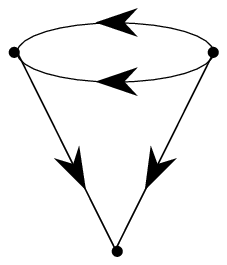}
\put(-1.3,2.7){$e_1$}
\put(-1.3,1.5){$e_2$}
\put(-2.5,1.1){$e_3$}
\put(-0.5,1.1){$e_4$}
\end{center}
\end{figure}

\noindent Then the involution is given by $\iota(e_1)=-e_2$ and
$\iota(e_3)=e_4$.
The elements
$$
h_1=e_1+e_3-e_4,\ h_2=e_2+e_3-e_4
$$
form a basis of $H_1(C,\ZZ)$ and with respect to this basis
$$
\iota=\begin{pmatrix} 0 & -1\\ -1 & 0\end{pmatrix}.
$$
It follows that
$$
X^-=\langle \frac 12 (h_1+h_2)\rangle;\  [X]^-=2X^-=\langle h_1+h_2\rangle.
$$
The picture below shows the Delaunay decomposition $\Delta$ of
$X_{\RR}$ and its
intersection with $X^{-}_{\RR}$.
%\newpage

\setlength{\unitlength}{1cm}
\begin{figure}[h]
\begin{center}
\epsfig{file=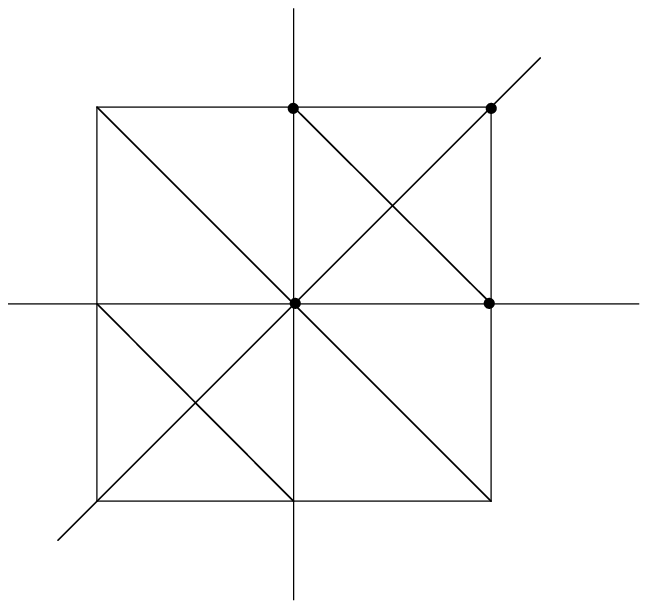}
\put(-3.4,5.2){$h_2$}
\put(-1.1,5.1){$X^-_{\RR}$}
\put(-1.3,2.7){$h_1$}
\end{center}
\end{figure}

\noindent The decomposition $[\Delta]^-$ is the Delaunay decomposition
with respect to
$X^-$. In particular it is {\em not} Delaunay with respect to $2X^-$. On
the other hand
condition (*) is fulfilled for $X^-$ since the rank of the lattice is 1.
Again in this
case $Y^-=2X^-$, since we can deform $(C,\iota)$ to a smooth curve
$(C_t,\iota)$ where the involution has no smooth fixed points and whose
Prym, therefore, has twice a principal polarization.
We want to use this example to discuss the relationship between
the degenerating families of the principally polarized and the twice
principally polarized Prym. The two families coincide outside the central
fiber, but have different central fibers, namely a nodal elliptic curve in
the first case and a union of $2\ \PP{^1}$'s intersecting in two points in
the second case. We start with a generic smoothing of the pair
$(C,\iota)$. In this case the mono\-dro\-my corresponds to the form
$B=\sum\limits^3_{j=1} z^2_j$. This defines a bilinear form on
$H_1(C, \ZZ)$ which, with respect to the basis $h_1, h_2$ of $H_1(C,
\ZZ)$, is given by the matrix
$$
B|_{H_1(C, \ZZ)}=\begin{pmatrix}3 & 2\\ 2 & 3\end{pmatrix}.
$$
Let $l=1/2(h_1+h_2)$ be a generator of $X^-$. The form
$[B/2]^-:Y^-\times X^-\rightarrow \ZZ$ is then given by
$[B/2]^-(l,2l)=5$. In order to describe degenerating families we have to
specify height functions $H_i:X^-\rightarrow \QQ; i=1,2$ which define the
decompositions $\Delta^-$ and $[\Delta]^-$. We define
$A^-:X^-\rightarrow\QQ$ by $A^-(ml)=5/4 m^2+1/2 m$ and set
$$
H_1(ml)=A^-(ml)+c_1(\bar{m})
$$
where
$$
c_1(\bar{m})=\left\{
\begin{array}{llcl}
0 &\operatorname{\quad if } m &\equiv & 0\mod 2\\
5/4 & \operatorname{\quad if } m &\equiv & 1 \mod 2,
\end{array}\right.
$$
respectively
$$
c_2(\bar{m})=\left\{
\begin{array}{llcl}
0 &\operatorname{\quad if } m &\equiv & 0\mod 2\\
1/4 & \operatorname{\quad if } m &\equiv & 1 \mod 2.
\end{array}\right.
$$
Note that only the quadratic part of $A^-$ is uniquely  determined by
$[B/2]^-$. Given
these data we can construct the degeneration families (for details
see(\cite{AN,A2}).
For this we consider the algebra $k[[t]][z][\zeta^{\pm}]$ where $t$ plays
the role of
the variable in the base and $z$ is a formal variable giving the grading.
We then
consider the subalgebras $\tilde{R_i}\subset k[[t]] [z] [\zeta^{\pm}]$
generated
by the elements $z\zeta^x t^{H_i(x)}, x \in X^-$. The group $Y^-$ acts on
$\tilde{R}_i$ by
$$
y^-:\ z\zeta^x t^{H_i(x)}\mapsto z\zeta^{x+y} t^{H_i(x+y)}.
$$
One then obtains the degenerating families as ${\cal A}_i=\Proj
(\tilde{R_i})/Y^-$. A
straight forward calculation shows that $\Proj (\tilde{R_i})$ is the torus
embedding
of a fan $\Delta_i$ which looks as follows
%\newpage
\setlength{\unitlength}{1cm}
\begin{figure}[h]
\begin{center}
\epsfig{file=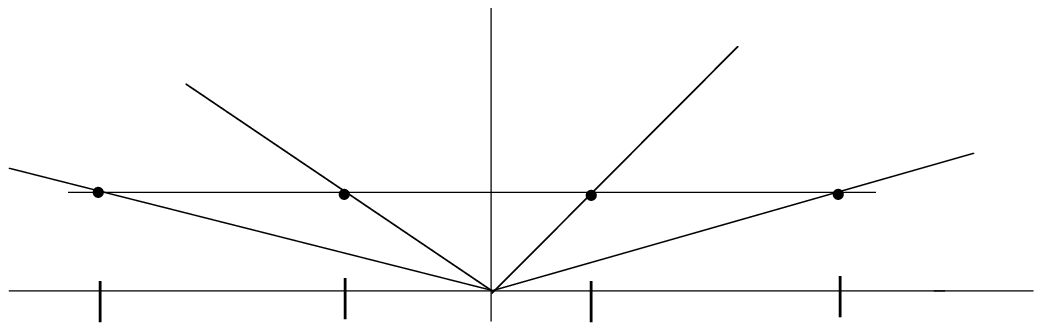}
\put(-5.4,2.5){$\RR_{\ge 0} f$}
\put(-4.6,-0.4){$2$}
\put(-2.0,-0.4){$7$}
\put(-7.2,-0.4){$-3$}
\put(-9.6,-0.4){$-8$}
\end{center}
\end{figure}

\noindent for $\Delta_1$, resp.
\newpage

\setlength{\unitlength}{1cm}
\begin{figure}[h]
\begin{center}
\epsfig{file=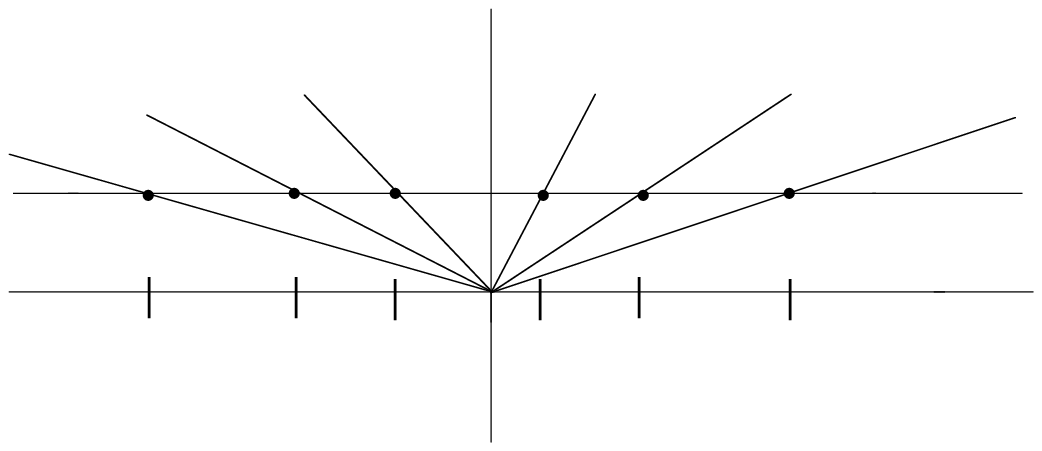}
\put(-5.3,4.0){$\RR_{\ge 0} f$}
\put(-5.1,0.9){$1$}
\put(-4.1,0.9){$3$}
\put(-2.6,0.9){$6$}
\put(-6.8,0.9){$-2$}
\put(-7.9,0.9){$-4$}
\put(-9.3,0.9){$-7$}
\end{center}
\end{figure}
\noindent for $\Delta_2$. The projection onto the base is given by
projecting onto
$\RR_{\ge 0} f$. {}From this we can read off that ${\cal A}_1$ has one
$A_5$-singularity
and that
${\cal A}_2$ has one $A_2$- and one $A_3$-singularity. The fan
\setlength{\unitlength}{1cm}
\begin{figure}[h]
\begin{center}
\epsfig{file=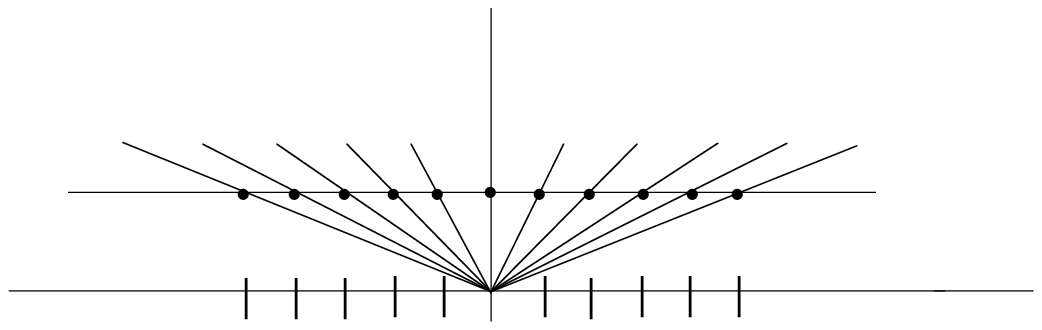}
\put(-5.6,-0.4){$0$}
\put(-5.1,-0.4){$1$}
\put(-4.6,-0.4){$2$}
\put(-4.1,-0.4){$3$}
\put(-3.6,-0.4){$4$}
\put(-3.1,-0.4){$5$}
\put(-6.3,-0.4){$-1$}
\put(-6.8,-0.4){$-2$}
\put(-7.3,-0.4){$-3$}
\put(-7.8,-0.4){$-4$}
\put(-8.3,-0.4){$-5$}
\end{center}
\end{figure}
%\newpage

\noindent is a common refinement and corresponds to blowing up ${\cal A}_1$ and
${\cal A}_2$ in
their singular points. The central fibre is then of type $I_6$, i.e. a
$6$-gon of
$\PP^{1}$'s.

%%%%%%%%%%%%%%
Note that all the examples which we have discussed so far lie in
$\overline{\cal R}_g$, provided we have no smooth fixed points.

\subsection{Mixed types of nodes}
\subsubsection{Branchwise fixed nodes and non-fixed nodes}
(a)\
Let $C$ be an irreducible curve with $2a$ non-fixed nodes $Q_1,\ldots, Q_{2a}$
and $b$ branchwise fixed nodes $Q_1,\ldots, Q_b$.
Again we choose the indexing such that
$\iota(Q_j)=Q_{2a-j+1}$. We
denote the edges corresponding to $Q_j$ by $e_j$ and those corresponding to
$Q_k$ by
$f_k$. Then $\iota(e_j)=e_{2a-j+1}$ and $\iota(f_k)=f_k$. Hence
$$
[H_1(\Gamma(\ZZ))]^-=[X]^-=\langle e_1-e_{2a},\ldots,e_{a}-e_{
a+1}\rangle=2X^-
$$
and the discussion is very similar to that of case (2.1). In particular
$[\Delta]^-$ is
Delaunay with respect to $2X^-$, conditions (**) and (*) are
fulfilled and $Y^-=2X^-$.

\noindent(b)\
We can also have a combination of  branchwise fixed and non-fixed nodes
if the
curve $C$ is
reducible. E.g., we can have $C=C_1\cup C_2$ where $C_1$ and $C_2$ intersect
in $2a$ nodes
$Q_1,\ldots, Q_{2a}$ which are non-fixed and $b$ branchwise fixed nodes
$Q_1,\ldots, Q_b$
which implies that necessarily $\iota(C_{i})=C_i$. For
simplicity we shall discuss the case $a=2$ and $b=1$. This implies the
existence of a
smooth fixed point on each of the components and we shall assume that there is
exactly one such point on each $C_i$. Then the graph
$\Gamma(C)$ of $C$ looks as follows

\setlength{\unitlength}{1cm}
\begin{figure}[h]
\begin{center}
\epsfig{file=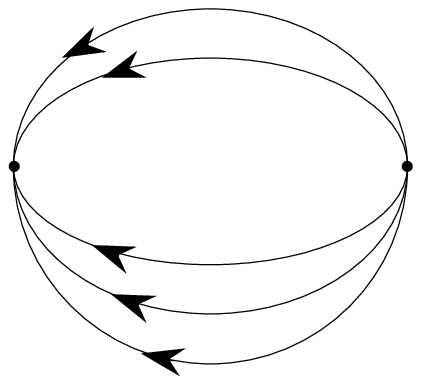,scale=0.8}
\put(-1.7,3.1){$e_1$}
\put(-1.7,2.3){$e_2$}
\put(-1.7,1.1){$e_3$}
\put(-1.7,0.65){$e_4$}
\put(-1.7,-0.3){$f$}
\end{center}
\end{figure}
%\newpage

\noindent where $f$ belongs to the branchwise fixed node $Q_1$.
As in the discussion
before we shall label
the nodes in such a way that $\iota(e_1)=e_4, \iota(e_2)=e_3$. Since $Q_1$
is branchwise fixed
we have $\iota(f)=f$. The elements
$$
h_1=e_1-e_3, h_2=e_4-e_2, h_3=f-e_1, h_4=f-e_4
$$
form a basis of $H_1(C,\ZZ)$ and with respect to this basis
$$
\iota=\begin{pmatrix}
\begin{pmatrix}
0 & 1\\ 1 & 0
\end{pmatrix} &\\
&
\begin{pmatrix}
0 & 1\\ 1 & 0
\end{pmatrix}
\end{pmatrix}.
$$
Let
$$
l_1=\frac 12(h_1-h_2),\  l_2=\frac 1 2 (h_3-h_4).
$$
Then
$$
[X]^-=\langle 2 l_1, 2 l_2\rangle = 2X^-.
$$
The dicing which appears in condition (*) now looks as follows

\newpage
\setlength{\unitlength}{1cm}
\begin{figure}[h]
\begin{center}
\epsfig{file=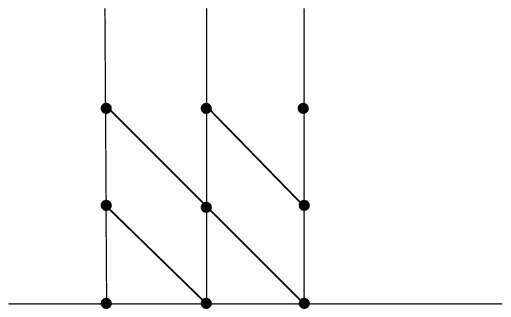}
\put(-4.9,2.2){$4l_2$}
\put(-4.9,1.2){$2l_2$}
\put(-3.3,-0.3){$2l_1$}
\put(-2.2,-0.3){$4l_1$}
\end{center}
\end{figure}
\noindent and this shows that (*) is fulfilled. We even have that
$[\Delta]^-$ is
Delaunay with
respect to $2X^-$, i.e. (**) holds.
Comparing this to example \ref{F-S} shows that the extra
node $Q_1$ ''improves'' the situation in an
essential way. The building blocks of the toric part are quadrics
$\PP^1\times\PP^1$.

\subsubsection{Swapping nodes and non-fixed nodes}
(a)\
We first assume that $C$ is irreducible and that it has $2a$ non-fixed nodes
$Q_1,\ldots,Q_{2a}$
and $b$ swapping nodes $R_1,\ldots, R_b$. We shall assume that
$\iota(Q_j)=Q_{2a-j+1}$
and correspondingly $\iota(e_j)=e_{2a+j-1}$. We will denote the edges
corresponding to the
nodes $R_j$ by $f_j$. The involution acts on these edges by
$\iota(f_j)=-f_j$. This case
is  related to case (\ref{Case 3.3}). We find
$$
\begin{array}{lcl}
[X]^- &=& \langle f_1,\ldots, f_b, e_1-e_{2a},\ldots, e_{a}-e_{a+1}\rangle,\\
2X^- &=& \langle 2f_1,\ldots,2f_b, e_1-e_{2a},\ldots,
e_{a}-e_{a+1}\rangle\neq[X]^- \mbox{ if } b\ge 1.
\end{array}
$$
The decomposition $[\Delta]^-$ is not Delaunay with respect to
$2X^-$.
Recall that swapping nodes contribute $2$ fixed points on a smoothing $C_t$.
We can compute the
polarization of $P(C_t, \iota)$ and compare this with Proposition \ref{prop18}.
We find that $|X^-/Y^-|=2$ if we have
no smooth fixed points and $X^-=Y^-$ if we have smooth fixed points.
In particular $Y^- \neq 2X^-$ unless $b=1$ and we have no smooth fixed
points. Strictly speaking we cannot speak about condition (*) since this
is not a point in $\oRg$, but the dicing condition still holds.
The building blocks of the toric part are all of the form
$(\PP^1)^b\times(\PP^1)^a$ and the number of these building blocks in
$[\overline{P}]$ is equal to $2^b$.\\
\noindent (b)\
Now assume that $C=C_1\cup C_2$. We shall assume that all nodes lie
on the intersection of $C_1$ and $C_2$.
If we want swapping nodes, then we
must necessarily
assume that $\iota(C_1)=C_2$ and this in turn rules out branchwise fixed nodes.
For the sake of simplicity we shall assume that we have no smooth
fixed points. Again we
assume that we have non-fixed nodes $Q_1,\ldots,Q_{2a}$ and swapping nodes
$R_1,\ldots,R_b$ .
If we choose the orientation in such a way that all edges start at
the same vertex, then
$\iota(e_j)=-e_{2a-j+1}$ and
$\iota(f_j)=-f_j$. To obtain a basis
of $H_1(C,\ZZ)$ we define
$$
\begin{array}{l}
h_0 = f_1-e_1\\
h_1 = e_{2a}-f_1\\
h_{2m} = e_m-e_{2a-m};\   m=1,\ldots, a-1\\
h_{2m+1} = e_{m+1}-e_{2a-m+1};\   m= 1,\ldots , a-1\\
w_k = f_l-f_1;  l=2,\ \ldots, b.
\end{array}
$$
With respect to this basis
$$
\iota=\begin{pmatrix}
\begin{pmatrix}
0 & 1\\ 1 & 0
\end{pmatrix} &&&&&\\
& \ddots &&&&\\
&&\begin{pmatrix}0 & 1\\ 1 & 0
\end{pmatrix}&&&\\
&&&-1&&\\
&&&&\ddots&\\
&&&&&-1
\end{pmatrix}
$$
where we have $a$ blocks of type $\begin{pmatrix}0 & 1\\ 1 & 0\end{pmatrix}$
and $b-1$ entries
$-1$. Let
$$
l_i=\frac 12 (h_{2i}-h_{2i+1});\ i=0,\ldots, a-1.
$$
Then
$$
%\begin{array}{l}
[X]^- = \langle 2 l_0,\ldots, 2 l_{a-1}, w_2, \ldots, w_b\rangle, \quad
X^- = \langle l_0,\ldots, l_{a-1}, w_2,\ldots, w_b\rangle.
%\end{array}
$$
The abelian part of the degenerate Prym carries twice a principal
polarization. On the other hand the Prym of a smoothing $(C_t,\iota)$
has a polarization of type $(1, \ldots,1,2, \ldots ,2)$
where the number of $1$'s is
$b-1$. Hence $|X^-/Y^-|=a$ and, in particular $X^- \neq 2Y^-$
if $b \geq 2$.

As a special case we consider $a=2, b=1$. Then the decomposition
$[\Delta]^-$
looks as follows
%\newpage
\setlength{\unitlength}{1cm}
\begin{figure}[h]
\begin{center}
\epsfig{file=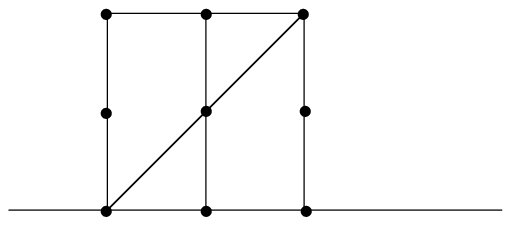}
\put(-4.7,1.8){$2l_1$}
\put(-4.7,0.9){$l_1$}
\put(-3.3,-0.4){$l_0$}
\put(-2.2,-0.4){$2l_0$}
\end{center}
\end{figure}

\noindent This shows that $[\Delta]^-$ is not Delaunay with respect to $2X^-$
The curve $(C,\iota)$ does not belong to a point in $\oRg$, but the dicing
condition of (*) is fulfilled.
The building blocks for the toric part of the
``middle'' Prym
$[\overline{P}]$  are $\PP^2$'s and
Hirzebruch surfaces $F_1$.

Finally we consider the case $a=1$ and $b=2$. In this case both $\Delta^-$ and
$[\Delta]^-$ looks as follows
\setlength{\unitlength}{1cm}
\begin{figure}[h]
\begin{center}
\epsfig{file=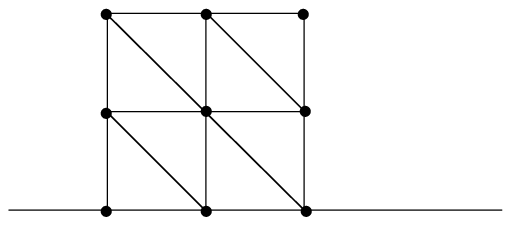}
\put(-4.7,1.8){$2w_2$}
\put(-4.7,0.9){$w_2$}
\put(-3.3,-0.4){$l_0$}
\put(-2.2,-0.4){$2l_0$}
\end{center}
\end{figure}

\noindent The dicing condition of (*) is fulfilled, but
$[\Delta]^-$ is clearly not Delaunay
with respect to $2X^-=\langle 2l_0, 2w_2\rangle$. The building blocks of the
toric part are projective planes.

\subsubsection{Branchwise fixed nodes and swapping nodes}
(a)
Assume that $C$ is irreducible and that it has $a$ branchwise fixed nodes
$Q_1,\ldots, Q_a$ and $b$ swapping nodes
$R_1,\ldots, R_b$.
Moreover we assume that there are $2r$ smooth fixed points.
Let $e_1,\ldots, e_a, f_1,\ldots f_b$ denote the
corresponding edges. Then $\iota(e_i)=e_i$ and $\iota(f_j)=-f_j$. Hence
$$
[X]^-=X^-=\langle f_1, \ldots, f_b\rangle.
$$
The decompositions $\Delta^-$ and $[\Delta]^-$ are the standard cubes with
vertices in $X^-$. In particular $\Delta^-$ is Delaunay with respect to
$X^-$ and the dicing
condition of (*) is fulfilled. Computing the degree of the polarization by
looking at a smoothing
we find that $|X^-/Y^-|=2$ if $r=0$ and $X^-=Y^-$ otherwise.

\noindent (b)\
Finally assume that $C=C_1\cup C_2$ and that all nodes lie in the
intersection of the irreducible components $C_1$ and $C_2$.
Then we cannot have branchwise fixed nodes and swapping nodes
simultaneously. Let us
assume that we have $b$ swapping nodes $Q_1,\ldots, Q_b$ with
corresponding edges
$f_1,\ldots, f_b$. Then $\iota(f_i)=-f_i$ and
$H_1(\Gamma,\ZZ)=[H_1(\Gamma,\ZZ)]^-$.
Again $[\Delta]^-$ is the Delaunay decomposition into standard cubes with
respect to
$[H_1(\Gamma,\ZZ)]^-=[X]^-$. As before we have $|X^-/Y^-|=2$ if we have no
smooth fixed points and $X^-=Y^-$ otherwise. In particular $Y^- \neq 2X^-$
if $b \geq 2$ or $b=1$ and we have smooth fixed points.

%%%%%%%%%%%%%%%%%%%%%%%%%%%%%%%%%%%%%%%%%%%%%%%%%%%%%%%%%%%%%%%%%%%%%%%%%%%%%%%
% References
%%%%%%%%%%%%%%%%%%%%%%%%%%%%%%%%%%%%%%%%%%%%%%%%%%%%%%%%%%%%%%%%%%%%%%%%%%%%%%%
\bibliographystyle{alpha}

\noindent
V.~Alexeev, valery@math.uga.edu\\
Department of Mathematics, University of Georgia\\
Athens, Georgia 30602, USA\\

\noindent CH.~Birkenhake, birken@mi.uni-erlangen.de\\
Fachbereich Mathematik, Universit\"at Mainz\\
Staudingerweg 9, D 55099 Mainz, Germany\\

\noindent K.~Hulek, hulek@math.uni-hannover.de\\
Institut f\"ur Mathematik, Universit\"at Hannover\\
Welfengarten 1, D 30060 Hannover, Germany\\

\end{document}